\documentclass[11pt]{amsart}
\usepackage{graphicx}
\usepackage{amssymb}
\usepackage{epstopdf}
\usepackage{amsmath}
\usepackage{amscd}
\usepackage{amsthm,amsfonts,amsxtra,amssymb,amscd}
\usepackage[all]{xy}
\usepackage{enumerate}

\usepackage{hyperref}

\theoremstyle{plain}
\newtheorem*{lemma*}{Lemma}
\newtheorem{lemma}[subsection]{Lemma}
\newtheorem*{theorem*}{Theorem}
\newtheorem{theorem}[subsection]{Theorem}
\newtheorem*{proposition*}{Proposition}
\newtheorem{proposition}[subsection]{Proposition}
\newtheorem*{corollary*}{Corollary}
\newtheorem{corollary}[subsection]{Corollary}
\theoremstyle{definition}
\newtheorem*{definition*}{Definition}
\newtheorem{definition}[subsection]{Definition}
\newtheorem*{example*}{Example}
\newtheorem{example}[subsection]{Example}
\theoremstyle{remark}
\newtheorem*{remark*}{Remark}
\newtheorem{remark}[subsection]{Remark}

\title[Index theory ]{Vector partition functions and index   of transversally elliptic operators.}
\author{
C. De Concini,\quad
C. Procesi,\quad M. Vergne}

\begin{document}

\begin{abstract}
Let $G$ be a torus acting linearly on a complex vector space $M$, and let $X$ be the list of weights of $G$ in $M$.
We determine the equivariant $K$-theory of the open subset $M^f$ of $M$ consisting of points with
finite stabilizers.  We identify it to the space $DM(X)$  of
 functions on the lattice $\hat G$,
  satisfying  the cocircuit difference equations associated to $X$,  introduced by Dahmen--Micchelli
   in the context of the theory of splines  in order to study vector partition functions
    (cf. \cite{DM}).

This  allows us to determine the range of the index map from
$G$-transversally elliptic operators on $M$ to generalized functions on $G$
and to prove that the index map is an isomorphism on the image.
This is   a setting studied by Atiyah-Singer \cite{At}  which is in a sense {\it universal}  for index computations.

\end{abstract}\maketitle

\section{Introduction}

\smallskip

In 1968 appears the fundamental work of Atiyah and Singer on the index theorem of elliptic operators, a theorem formulated in successive steps of generality (\cite{Ats}, \cite{Ats2}). One general and useful setting is for operators  on a manifold $M$ which satisfy a symmetry with respect to a compact Lie group  $G$  and are elliptic in directions transverse to the $G$-orbits.
The values of the index are generalized functions on $G$.

In his Lecture Notes \cite{At} describing joint work with I.M. Singer,
Atiyah explains how to reduce general computations to the case in which $G$ is a torus,
 and the manifold $M$ is a complex linear representation $M_X=\oplus_{a\in X}L_a$, where  $X\subset \hat G$ is a finite list of characters and  $L_a$  the one dimensional complex line where $G$ acts by the character $a\in X$.  He then computes explicitly in several cases and ends his introduction saying

"  ... for a circle (with any action)  the results are also quite explicit.  However  for the general case we give only a
reduction  process and one might  hope for something explicit.  This probably requires  the development of an appropriate
algebraic machinery, involving cohomology but    going beyond it."

The purpose of this paper is to provide  this  algebraic machinery,  which turns out to be a spinoff of the theory of splines, and complete this computation.
In fact   we show that, once one has guessed  the
right formulas, the actual
computations are obtained using directly the basic properties of $K$-
theory and turn out to be quite simple. In particular they neither
involve computations of the index of specific operators nor
any sophisticated commutative algebra.
Let us now explain in some detail  the content  of the paper.\smallskip

We denote by $T^*_GM$ the closed subset of $T^*M$, union of the conormals to the $G$ orbits.
The symbol $\sigma(x,\xi)$ of a pseudo-differential  transversally elliptic operator $A$  on $M$ determines
an element of the equivariant $K-$theory group  $K_G^0(T^*_GM)$  which is defined in a  topological fashion.
The analytic index of the operator $A$ is  the virtual trace class representation of $G$  obtained as difference  of the spaces of
solutions  of $A$ and its adjoint $A^*$ in an appropriate Sobolev  space.
This index depends only of the class defined by $\sigma$ in $K_G^0(T_G^*M)$, so that
the  index defines a $R(G)$ module homomorphism from $K^0_G(T^*_GM)$ to virtual trace class representations of $G$.

Let $V$ be the dual of the Lie algebra of $G$, so that $\hat G$ is identified with the set $\Gamma$ of weights of $G$, a lattice in $V$. Denote by $\mathcal C[\Gamma]$ the space of $\mathbb Z$ valued functions on $\Gamma$.
A trace class representation  $\Theta$ of $G$ can be decomposed as a  direct sum
$\Theta=\oplus_{\gamma\in \Gamma}f(\gamma)L_ \gamma$
, where $f\in\mathcal C[\Gamma]$ is a function on $\Gamma$ with at most polynomial growth. We say that $f$ is the multiplicity function of $\Theta$.
We also denote by $\Theta(g)$ the trace of the representation $\Theta$. It is a generalized function on $G$.
In this article, we determine the  subspace $\tilde{\mathcal F}(X)$ of $\mathcal C[\Gamma]$
arising
from multiplicities of indices of transversally elliptic operators on $M_X$.
In fact, we  study the more general case of a  non necessarily connected compact abelian group $G$, as,
in recurrence steps, non connected subgroups of $G$  appear.
However  in this introduction, we stay with the case of a torus, that is a direct product of circle groups $S^1$.
It is also harmless to assume that $G$ acts on $M$ without non zero fixed vectors, and we do so.

In  \cite{At}, the image of the index map  is  described explicitly  for $S^1$
and implicitly in general by a constructive algorithm.
It is useful to recall Atiyah-Singer's  result for the simple  case where
$G=S^1$ acts by homotheties on $M_X=\mathbb C^{k+1}$.
We  denote by $t$ the basic character of $S^1:=\{t\,|\,\,|t|=1\}$, so that $R(G)=\mathbb Z[t,t^{-1}]$  and $X=[t,t,\ldots,t]$ , $k+1$ times.

First Atiyah-Singer  constructed  a ``pushed" $\overline \partial$ operator on $M_X$, with index the trace of the
representation of $G$ in the symmetric algebra $S(M_X)$.
As
$\binom{n+k}{k}=\frac{(n+1)\cdots (n+k)}{k!}$ is the dimension of the  space of homogeneous polynomials  in $k+1-$variables and degree $n$,
the character of $S^1$ in $S(M_X)=\oplus_{n=0}^{\infty} S^n(M_X)$
is the generalized function
$$ \Theta_X(t):=\sum_{n=0}^\infty  \binom{n+k}{k}t^n. $$

Remark that $n\mapsto  \binom{n+k}{k}$ extends to a polynomial function on $\mathbb Z$.
For any $n$ positive or negative,
the  function $n\mapsto  \binom{n+k}{k}$  represents the dimension of a {\it virtual space}, the
alternate sum of the  cohomology spaces  of the sheaf  $\mathcal O(n)$ on $k-$ dimensional projective space.
In particular, the  tangential Cauchy-Riemann operator on the unit sphere $S_{2k+1}$ of $\mathbb C^{k+1}$  is a transversally elliptic operator with index
$$\theta_X(t):=\sum_{n=-\infty}^\infty  \binom{n+k}{k}t^{n}, $$
a generalized function on $G$ supported at $t=1$.
Then it is proved in \cite{At} that the index map is an isomorphism  from $K_G^0(T^*_G M_X)$ to the space $\tilde {\mathcal F}(X)$ of generalized functions on $G$ generated by $\Theta_X$ and $\theta_X$ under multiplication by an element of $R(G)=\mathbb Z[t,t^{-1}]$.
In fact the $R(G)$ module generated by $\Theta_X$ is free over $R(G)$.
The $R(G)$ module generated by
$\theta_X$ is the torsion submodule. This submodule is the module of polynomial functions on $\mathbb Z$ of degree at most $k$ so it is a free $\mathbb Z$-module of rank $k+1$.  It corresponds to
indices of operators on $\mathbb C^{k+1}-\{0\}$, the set where $S^1$ acts freely and is the space of solutions of the difference equation $\nabla^kf=0$ where we define the operator $\nabla$ by $(\nabla f)(n)=f(n)-f(n-1)$.

In higher dimensions, the single difference equation $\nabla^kf=0$  must be replaced by a system of difference equations, discovered by Dahmen--Micchelli as the natural generalization of a system of differential equations  associated to splines  in approximation theory.

This is
the system of  difference equations   $\nabla_Yf=0$ associated to cocircuits $Y$ of $X$ ($X\subset \Gamma$,  a sequence   of weights of $G$).

Let us   recall the definitions introduced in \cite{dpv}.
\begin{itemize}
\item A
subspace $\underline r$ of $V$ is called {\it rational} if $\underline r$ is
the span of a sublist of $X$.
\item A  {\it
cocircuit} $Y$ in $X$  is a sublist  of the form $Y=X\setminus H$  where $H$ is a rational
hyperplane.
\item Given $a\in\Gamma$, the {\it difference operator}
$\nabla_a$ is the operator on functions $f$ on $\Gamma$ defined by
$\nabla_a(f)(b):=f(b)-f(b-a)$.\item  For a list $Y$ of vectors, we set
$\nabla_Y:=\prod_{a\in
Y}\nabla_a$.
\item The Dahmen-Micchelli space introduced in \cite{DM} is the space
$$DM(X):=\{f\,|\,  \nabla_Y
f=0,\ \\ \text{ for every  cocircuit $Y$ in }X\}.$$
\end{itemize} It is easy to
see that $DM(X)$ is finite dimensional and consists
of quasi--polynomial
functions. We describe it in detail in \S \ref{TarPa}.  In  \cite{dp1},
we introduced non homogeneous difference equations   and  defined
$$\mathcal F(X):=\{f\,|\,  \nabla_{X\setminus  \underline r }
f \text{ is supported on }  \underline r \text{ for every
rational subspace    } \underline r  \}.
$$
Clearly  $DM(X)$ is contained in
$\mathcal F(X)$.

We denote by $\tilde{\mathcal F}(X)$  the space of functions
 generated by ${\mathcal F}(X)$ under translations by elements of $\Gamma$. Thus,
$\tilde{\mathcal F}(X)$  is the space of functions  $f\in \mathcal C[\Gamma]$  such that $\nabla_{X\setminus  \underline r }
f$  is supported on a finite number of $\Gamma$ translates of  $ \underline r$  for every
proper rational subspace    $ \underline r$.
Both $\mathcal F(X)$
and $\tilde {\mathcal F}(X)$ have a very precise description given by Formulas  (\ref{decF}) and  (\ref{decFt}).
Our main theorem is:

\begin{theorem}\label{main}
Let $M_X:=\oplus_{a\in X} L_a$. Then the index multiplicity  map induces an
isomorphism from $K^0_G(T^*_G M_X)$ to $\tilde {\mathcal F}(X)$.
\end{theorem}

Our proof is inspired by Atiyah-Singer  strategy. We start by  determining  the equivariant $K$-theory
of
 $T_GM_X^f$,   $M_X^f$ being the open subset of $M_X$ of points with finite stabilizers. Here we assume that $X$ spans $V$. If  $s$ denotes the dimension of $G$, we prove

 \begin{theorem}

\begin{enumerate}
\item
We have isomorphisms

$$K^s_G(  M_X^f ) \cong K^0_G(T_G ^*M_X^f ) \stackrel{ind_m}\cong  DM(X ) $$
where $ind_m$ is the index multiplicity map.
\item
$K^{s+1}_G(M_X^f)=0$.
\end{enumerate}
\end{theorem}

We then study the equivariant $K$ -theory of the space $T_G^*M_X$ by induction over the natural stratification of $M_X=\cup_iF_i$, using as strata $F_i$ the union of all orbits of dimension $i$.

Let us comment on our method.
The proof of the exactness of several important sequences is
based on the purely combinatorial fact that, if $X:=[a,Z]$ is formed by adding a vector $a$ to a list $Z$, then the map
$$\nabla_a: DM(X)\to DM(Z)$$ is surjective and its kernel is also a space $DM(\tilde Z)$.  This is just one of the several incarnations of the deletion and restriction principle in the theory of hyperplane arrangements (see \cite{zas}).  It also corresponds to a natural geometric description of the zonotope  $B(X)$ (cf. \cite{dp1})  which we do not discuss here.

We give
an elementary proof of this result  based on a construction of generators for $DM(X)$.
This is done in the first two sections.
Although hidden in the proof, this system of generators
 corresponds to the generators given by Atiyah-Singer of $K_G^*(M_X^f)$.

We can deduce easily     generators for $K_G^0(T^*_G M_X)$ from Theorem \ref{main}.
Consider all possible complex structures $J$  on $M_X$ obtained in the following way. We consider
any decomposition of $X$ in $A$ and $B$  so that the cone $C(F):=C(A,-B)$ spanned by $A$ and $-B$  is a pointed cone. We then take the conjugate complex structures on the lines $L_b$, with $b\in B$.
Then, to this complex structure $J$ is associated a ``pushed" transversally elliptic $\overline \partial$
operator on $M_J$  with symbol
$At_J\in K_G^0(T^*_G M_X)$.
We recall the definition of $At_J$ in (\ref{At}).
As a consequence of our theorem, we obtain a simple proof of
Atiyah-Singer  system of generators for $K^0_G(T^*_G M_X)$.
\begin{theorem}
The symbols $At_J$ generates  the $R(G)$ module $K_G^0(T^*_GM_X)$.
\end{theorem}
The generalized functions  associated to $DM(X)$ are  distributions supported on a certain finite subset $P(X)$ of $G$, namely  the points $p\in G$ having a fixed point in $M_X^f$.
We will pursue in a subsequent article the description of the equivariant cohomology of $M_X^f$ and of the natural strata of $M_X$,  and relate this to the local components of the distributions associated to $DM(X)$ at a point $p\in P(X)$ and other fixed points.

\section{Combinatorics}

\subsection{Notations\label{Not}}

Let $G$  be a compact abelian Lie group of dimension $s$  with character group $\Gamma$.
If $\gamma\in \Gamma$, we denote by $\chi_\gamma(g):= g^\gamma$ the corresponding homomorphism $\chi_\gamma:G\to S^1$, where  $S^1:=\{u\in \mathbb C\,|\,\,|u|=1\}$.

This suggestive notation  for characters is justified by the following. If $H=(S^1)^N$ is the product of $N$ copies of the circle,  then every character of $H$ is of the form $(u_1,u_2,\ldots, u_N)\to u_1^{\gamma_1}\cdots u_N^{\gamma_N}$  where the $\gamma_i$ are integers. We denote this character by the multi-index notation  $u\mapsto u^\gamma$, where $\gamma=[\gamma_1,\gamma_2,\ldots, \gamma_N]\in \mathbb Z^N$.
Take now a closed subgroup $G$ of $H$, then every character of $G$ is also of this form where now  $\gamma$ is in a quotient $\Gamma$ of the lattice $\mathbb Z^N$.

We consider here $G$ as an abstract compact abelian Lie group, so that $\Gamma$ is a discrete finitely generated abelian group.   The group law  on  the character group $\Gamma$ of $G$  will be denoted additively.
If $\gamma_1,\gamma_2\in \Gamma$, the character associated to $\gamma_1+\gamma_2$ is the pointwise  multiplication $$\chi_{\gamma_1+\gamma_2}(g)=\chi_{\gamma_1}(g)
\chi_{\gamma_2}(g)=g^{\gamma_1} g^{\gamma_2}=g^{\gamma_1+\gamma_2}$$ of characters.

Let $U:=\mathfrak g$  be the Lie algebra of $G$ and    $V:=\mathfrak g^*$ the dual vector space. The spaces $U,V$ are real vector spaces of dimension $s$.
Denote by  $G^0$   the connected component of the identity   of $G$, so that $G^0$ is a connected torus with Lie algebra $U$, the group $G/G^0$ is a finite group,  and $G$ is isomorphic, although not in a canonical way, to $G^0\times G/G^0$.
We denote   by $\Gamma_t$ the torsion subgroup of $\Gamma$, identified with the character group of the group  $G/G^0$.  Then  the quotient group $\Gamma/\Gamma_t$  may  be identified  with the weight lattice $\Lambda$ of $G^0$, a lattice in $V$,   as follows.
Let $\gamma\in  \Gamma$ and  $g\in G^0$ and  write $g=\exp( X)$, where $X\in \mathfrak  g$. Then we can write
  $g^{\gamma}=e^{2i\pi \langle \lambda,X\rangle }$, where $\lambda\in \Lambda$ is an element of the weight lattice.
  In other words,  the linear form $\lambda\in \Lambda$ is the differential $d\chi_\gamma: \mathfrak g\to {\rm Lie}(S^1)\sim \mathbb R$  of the character
  $\chi_\gamma: G\to S^1$, where the exponential map $\mathbb R\to S^1$ is given  by  $u\to e^{2i\pi u}$.

If $\gamma \in \Gamma$, then  we denote by  $\overline \gamma\in V $
the corresponding element of the weight lattice $\Gamma/\Gamma_t\equiv \Lambda \subset V$.

Let  $\mathcal C[\Gamma]$  be  the space of functions  on $\Gamma$ with values in $\mathbb Z$. We consider $\mathcal C[\Gamma]$ as a $\Gamma$ module, by  associating to $\gamma\in \Gamma$ the translation operator $\tau_\gamma f(\nu):=f(\nu-\gamma)$. We set $\nabla_\gamma=1-\tau_\gamma$  and call it a {\it difference operator}.

Given a finite list $Y$ of  elements of $\Gamma$,  we set   $\nabla_Y:=\prod_{a\in Y} \nabla_a$. Also we set
  $a_Y:=\sum_{a\in Y}  a$ and
$\overline a_Y:=\sum_{a\in Y}  \overline a$.

\section{The space $DM(X)$\label{TarPa}}
We introduce here the basic space of functions   $DM(X)$  on $\Gamma$, which can be thought of as a multi dimensional generalization of suitable binomial coefficients.\bigskip

  Let us fix a finite non empty list $X$ of   characters in $\Gamma$.   Associated to $X$, we get  a list of  vectors $\overline X:=\{\overline a\,|\,a\in X\}$ in $\Lambda\subset V$.

Any non--zero element $\overline a\in \overline X$ defines a hyperplane in $U$. Thus,  we associate to $X$   the hyperplane  arrangement $\mathcal H_X$  in $U$  whose elements are all the intersections of the hyperplanes associated to the  non--zero elements of $\overline X$.
Given a set or list of vectors $A$ in $V$, we denote by $\langle A\rangle $  the   subspace spanned  by $A$.

\begin{definition}A subspace  $\underline r$ of $V$ is called {\it rational} (relative to $X$) if it is spanned by a sublist $\overline Y$ of $\overline X$.

The set of all rational subspaces  associated to $X$ will be denoted by $S_X$,  and
the subset  of rational subspaces of
dimension $i$ by $S_X^{(i)}$.

\end{definition}

If $Y$ is a sublist in  $X$,  the rational subspace spanned  by $\overline Y$ will be sometimes denoted, by abuse of notations, by  $\langle   Y\rangle:=\langle \overline Y\rangle $.
Similarly, by abuse of notations, we will often transfer to $X$ definitions  which make  sense   only for $\overline X$. For example,  if   $a\in X$  and $F$ is a    subset of $U$, we will say that {\it $a$ is positive on $F$} if $\overline a$ is positive on  every element $u$ of $F$.  If $A$ is a list of elements of $X$, we will say that {\it $A$ spans a pointed cone} if $\overline A$ spans a pointed cone. We then denote by $C(A)$ the elements of $\Gamma$ which are non negative integral linear combinations of elements of  $A$.
If $\underline r$ is a rational subspace of $V$, we write $X\cap \underline r$ (resp. $X\setminus \underline r$)  for the set of elements $a\in X$ such that $\overline a\in \underline r$ (resp.  $\overline a\notin \underline r$).

\begin{definition} If $  \underline r$ has codimension 1,  that is  $\underline r$  is  a rational hyperplane,   the sublist $X\setminus \underline r$ will be  called a {\it cocircuit}.
\end{definition}

Our list $X$   of elements in $\Gamma$  defines a set of subgroups of $G$. If $a\in X$, we denote by $G_a$ the kernel of the corresponding character $g\mapsto  g^a: G\to S^1$.
  The  set of all connected components of all intersections of these subgroups  is a finite family of cosets of certain tori and is called the {\it real toric arrangement} associated to $X$. We let $\mathcal T_X^i$  be the set of  tori of  dimension $i$ contained in this arrangement.

Notice that there is a 1--1 correspondence between the  tori in  $\mathcal T_X^i$ and the  rational subspaces of codimension $i$ by associating to $T$ the space $\langle X_T\rangle$, $X_T$  being the characters in $X$ containing $T$ in their kernel.

We are ready to
introduce the space of functions $DM(X)$.  This space is a straightforward  generalization of the space introduced by Dahmen--Micchelli for vectors in a lattice.   The extension to any finitely generated abelian  group is very convenient for   inductive  constructions as we shall see in Theorem \ref{LABAS}.

\begin{definition}
The space $DM(X)$  is the space of functions  $f\in \mathcal C[\Gamma]$  which are solutions of the linear equations  $\nabla_{X\setminus \underline r}f=0$ as $\underline r$ runs over all rational hyperplanes  contained in $V$.
\end{definition}
In other words,
the space $DM(X)$  is the space of functions  $f\in \mathcal C[\Gamma]$  which are solutions of the linear equations  $\nabla_{Y}f=0$ as $Y$
runs over  sublists of $X$ such that $X\setminus Y$ do not generate $V$

In particular, if $X'$ is a sublist of $X$, the space $DM(X')$ is a subspace of $DM(X)$.

\begin{example} If   $G=S^1\times \mathbb Z/2\mathbb Z$, $\Gamma=\mathbb Z\times  \mathbb Z/2\mathbb Z$. Take the list
$X=\{1+\epsilon,2+\epsilon\}$, where $n=1,2$  denotes  the character $\chi_n(t)= t^n$ of $S^1$  and $\epsilon$ the non trivial character of $\mathbb Z/2\mathbb Z$. Then $DM(X)$  consists of the functions $f$ such that for every $\gamma\in\Gamma$
$$f(\gamma)-f(\gamma-1-\epsilon)-f(\gamma-2-\epsilon)+f(\gamma-3)=0.$$
It is then clear that $f$ is uniquely determined by the values it takes on the points $0,-1,-2,-\epsilon,-1-\epsilon,-2-\epsilon$ and it is a free abelian group of rank 6.
	\end{example}

Notice that, if $\Gamma=\Gamma_t$ is finite, then $DM(X)=\mathbb Z[\Gamma]$ and does not depend on $X$. Also,  if     the joint kernel of all the elements in $X$ is of positive dimension,  $DM(X)=0$ since $\nabla_{\emptyset}=\nabla_{X\setminus\langle X\rangle}=1$.
In view of this, we are often going to tacitly make the {\it non degeneracy} assumption that   $\langle X\rangle= V\neq \{0\}$.

Under this assumption, $DM(X)$ is the space  of solutions of the equations
$\nabla_Yf=0$ as $Y$ runs over all cocircuits.

\begin{remark}\label{real}
If $X=[A,B]$ is decomposed into two lists, and if $X'=[A,-B]$ is obtained from the list $X$
 by changing the signs of the elements $b$ in the sublist $B$, then  $DM(X)=DM(X')$. \end{remark}

\subsection{Convolutions in the space $\mathcal C[\Gamma]$ }The purpose of this and the following section is to exhibit explicit elements in $DM(X)$, which have a natural interpretation in index theory.

If $\gamma\in \Gamma$,  we denote by $\delta_\gamma$ the function on $\Gamma$ identically equal
to $0$ on $\Gamma$,  except  for $\delta_\gamma(\gamma)=1$.
With this notation, we also write an element $f\in  \mathcal C[\Gamma]$ as
$$f=\sum_{\gamma\in \Gamma}f(\gamma)\delta_\gamma.$$

The {\it support } of  a function    $f\in \mathcal C[\Gamma]$  is the set of elements  $ \gamma\in \Gamma$ with $f(\gamma)\neq 0$.
If $S$ is a subset of $V$, we will say, by abuse of notations,  that  $f$ is supported on
$S$ if $f(\gamma)=0$ except if $\overline \gamma\in  S$.

\begin{remark}\label{supconv}
If $S_1,S_2$ are the supports of $f_1,f_2$, the
convolution product  $f_1*f_2$ is defined
  when,  for every  $\gamma\in \Gamma$,  we have only finitely many pairs $(\gamma_1,\gamma_2)\in S_1\times S_2$ with  $\gamma=\gamma_1+\gamma_2$.
 Then  $$(f_1*f_2)(\gamma)=\sum_{\gamma_1\in S_1,\gamma_2\in S_2\,|\, \gamma_1+\gamma_2= \gamma }f_1(\gamma_1)f_2(\gamma_2).$$
\end{remark} In particular,  if $\gamma \in \Gamma$,  $(\delta_{\gamma}*f)(x)=f(x-\gamma)=(\tau_{\gamma}f)(x)$. Thus convolution with  elements of finite support,  i.e. elements in  $\mathbb  Z[\Gamma]=R(G)$, gives the $R(G)$ module structure on  $\mathcal C[\Gamma]$
already introduced.

Recall that,  for $\gamma\in \Gamma$, we denoted by $g\to g^\gamma$ the corresponding function on $G$.

If  $f\in \mathcal C[\Gamma]$  has finite support, the function  $g\mapsto \sum_{\gamma\in \Gamma}f(\gamma)g^{\gamma}$
is a complex-valued function on $G$ with Fourier coefficients $f(\gamma)$.
When $f$ is arbitrary, we also write the formal series of functions on $G$ defined by
$$\Theta(f)(g)=\sum_{\gamma\in \Gamma} f(\gamma) g^{\gamma}.$$
This is usually formal although, in the cases of interest in this paper,  it may give a convergent series in the sense of generalized functions on $G$.
We denote by $C^{-\infty}(G)$ the space of generalized
functions on $G$.  The space $C^{\infty}(G)$ of
smooth functions on $G$ is naturally a subspace of
$C^{-\infty}(G)$. We will often use the notation  $\Theta(g)$ to
denote a generalized function $\Theta$ on $G$ although (in
general) the value of $\Theta$ on a particular point $g$ of $G$
does not have a meaning. By definition, $\Theta$ is a linear form
on the space of smooth densities on $G$ (satisfying suitable continuity conditions).

In fact, all the elements  $f=\sum_{\gamma}f(\gamma)\delta_\gamma$ that we shall consider have the property that the function  $f(\gamma)$ has polynomial growth at infinity. This implies that  the series $\Theta(f)$ converges in the distributional sense, and $f$ can be interpreted   analytically as the Fourier series of  a $C^{-\infty}$ function  $\Theta(f)$ on $G$.
The convolution product on $f_1,f_2$ (if defined) corresponds to the multiplication (if defined) of
the functions $\Theta(f_1)(g)\Theta(f_2)(g)$.

If $\gamma\in \Gamma$ is of infinite order, we set:
$$H_\gamma :=\sum_{k=0}^\infty \delta_{k\gamma}.$$
The function $H_\gamma$ is supported on the ``half line"  $\mathbb Z_{\geq 0}\gamma$
and is the discrete analogue of the Heaviside function.  In fact  $\nabla_\gamma H_\gamma=\delta_0.$\smallskip

Let  $A=[a_1,a_2,\ldots, a_n]$  be a finite list of vectors in $\Gamma\setminus \Gamma_t$ such that
$\overline A$ spans a pointed cone in $V$, and consider the {\it generalized cone}
$$C(A):=\{\sum_{i=1}^n m_ia_i,\ m_i\in\mathbb Z_{\geq 0}\}\subset\Gamma.$$
which, according to our conventions, we treat as {\it  a pointed cone}.

We easily see that every element $\gamma\in C(A)$  can be written in the form $\gamma=\sum_{i=1}^n m_ia_i,\ m_i\in\mathbb Z_{\geq 0}$ only in finitely many ways. It follows that the convolution $H_A$ given by
$$H_A:=H_{a_1}*H_{a_2}* \cdots* H_{a_n}$$
    is well defined and supported on  $C(A).$

\subsection{The   function ${\mathcal P}_{X\setminus \underline r}^{F_{\underline r}}$\label{convoli}}
Given a rational subspace $\underline r$,    $X\setminus \underline r$ defines a hyperplane arrangement in the space   $\underline r^\perp\subset U$ orthogonal to $\underline r$. Take an open  face $F_{\underline r}$ in   $\underline r^\perp$ with respect to this hyperplane arrangement.   By definition a vector  $u\in \underline r^\perp$ and such that $\langle u,\overline a\rangle \neq 0$ for all $a \in X\setminus \underline r$
lies in a unique such face $F_{\underline r}$.

 We get a decomposition of $X\setminus \underline r$ into the two lists $A,B$ of elements which are positive (respectively negative) on $F_{\underline r}$.  We denote   by $C(F_{\underline r},X)$ the cone $C(A,-B)$
generated by the list  $[A,-B]$.
Then $C(A,-B)$ is a pointed cone.

Notice that, given $a\in\Gamma$, we have  $\nabla_{ a}=-\tau_{ a}\nabla_{-a}$.

We are going to consider the function  ${\mathcal P}_X^{F_{\underline r}}$
which is characterized by the following two properties.

\begin{lemma}\label{lafunteta}
There exists a unique element  ${\mathcal P}_{X\setminus \underline r}^{F_{\underline r}} \in \mathcal C[\Gamma]$
such that \begin{enumerate}\item
$(\prod_{a \in X\setminus\underline r }
\nabla_a) {\mathcal P}_{X\setminus \underline r}^{F_{\underline r}}=\delta_0
$.
\item
${\mathcal P}_{X\setminus \underline r}^{F_{\underline r}}$ is supported in $-b_B+C(A,-B)$. \end{enumerate}
\end{lemma}
\begin{proof}
We define   ${\mathcal P}_{X\setminus \underline r}^{F_{\underline r}}$ by  convolution product of
``Heaviside functions":
\begin{equation}\label{esplicit}
{\mathcal P}_{X\setminus \underline r}^{F_{\underline r}}= (-1)^{|B|}\delta_{-b_B}* H_A* H_{-B}
\end{equation}
where $b_B=\sum_{b\in B}b$.
As $(\prod_{a \in X\setminus\underline r }
\nabla_a)=
(-1)^{|B|}\tau_{ b_B}(\prod_{a \in A }
\nabla_a) (\prod_{b \in B }
\nabla_{-b})$, we see
that   ${\mathcal P}_{X\setminus \underline r}^{F_{\underline r}}$  is well defined,
satisfies the two properties and is unique.
\end{proof}

\begin{remark}
When $G$ is connected so that $\Gamma$ is a lattice,
 this function was constructed in \cite{dpv}.
\end{remark}

It is easily seen that  the series $\sum_{\gamma\in \Gamma}{\mathcal P}_{X\setminus \underline r}^{F_{\underline r}}(\gamma) g^{\gamma}$ defines a generalized function
$\Theta_{X\setminus \underline r}^{F_{\underline r}}$ on $G$  such that
 $\prod_{a\in X\setminus \underline r}(1-g^a)
 \Theta_{X\setminus \underline r}^{F_{\underline r}}(g)=1$ on $G$.
Thus  $\Theta_{X\setminus \underline r}^{F_{\underline r}}(g)$ is
the inverse of   $\prod_{a\in X\setminus \underline r}(1-g^a)$,
in the space of  generalized functions on $G$, with Fourier coefficients  in the ``pointed cone" $C(A,-B)$.

We can generalize this as follows.
Choose an orientation of each  $\underline r\in S_X$.
Take a pair of rational subspaces  $\underline r \in   S_X^{(i)}$ and $\underline t  \in   S_X^{(i+1)}$ with   $\underline r \subset \underline t . $  We say
that a vector $v$ in $\underline t \setminus \underline r$ is positive if the orientation on $\underline t$ induced by $v$ and the orientation on $\underline r$
coincides with that  chosen of $\underline t$.
Set
$$A=\{a\in X\cap \underline t\,|\,\overline a\text {\ is positive}\},\quad B=\{b\in X\cap \underline t\,|\,\overline b \text {\ is negative}\}.$$\
We define
\begin{equation}\label{esplicit1}{\mathcal P}^{\underline t+}_{ \underline r} :=(-1)^{|B|}\delta_{-b_B}*H_A*H_{-B},\end{equation}
\begin{equation}\label{esplicit2}
{\mathcal P}^{\underline t-}_{ \underline r} :=(-1)^{|A|}  \delta_{-a_A}*H_B*H_{-A}.\end{equation}

The function   ${\mathcal P}^{\underline t+}_{ \underline r}$ is supported  on the   cone $-\sum_{b\in B}b+ C(A,-B)$  while
${\mathcal P}^{\underline t-}_{ \underline r}$ is supported  on  the cone $-\sum_{a\in A}a+C(-A, B)$.

\begin{definition}\label{defiQ}
Take a pair of  oriented rational subspaces  $\underline r \in   S_X^{(i)}$ and $\underline t  \in   S_X^{(i+1)}$ with   $\underline r \subset \underline t . $
We define $${\mathcal Q}^{\underline t}_{ \underline r}(X)= {\mathcal P}^{\underline t+}_{ \underline r}- {\mathcal P}^{\underline t -}_{ \underline r}.$$

If $X$ is fixed, we will  write simply  $ {\mathcal Q}^{\underline t}_{ \underline r}$ instead of
${\mathcal Q}^{\underline t}_{ \underline r}(X)$.
\end{definition}

Similarly   ${\mathcal P}^{\underline t\pm }_{ \underline r}, {\mathcal Q}^{\underline t}_{ \underline r}$
define generalized functions on $G$.
\begin{definition}
We set
\begin{equation}\label{esplicit4}
\theta^{\underline t\pm }_{ \underline r}(g) := \sum_{\gamma\in \Gamma} {\mathcal P}^{\underline t\pm}_{ \underline r}(\gamma) g^{\gamma},
\end{equation}

\begin{equation}\label{esplicit3}
\theta^{\underline t}_{ \underline r}(g) :=
\theta^{\underline t+}_{ \underline r}(g) -\theta^{\underline t-}_{ \underline r}(g)
=\sum_{\gamma\in \Gamma} {\mathcal Q}^{\underline t}_{ \underline r}(\gamma) g^{\gamma}.
\end{equation}
\end{definition}

Let us use the notation  $D_Y(g):=\prod_{a\in Y}(1-g^a)$.
Then
\begin{equation}\label{equ}
D_{(X\setminus \underline t)\setminus  \underline r}(g)
\theta^{\underline t\pm }_{ \underline r}(g)=1,
\end{equation}

so that
$$D_{(X\setminus \underline t)\setminus \underline r}(g)
\theta^{\underline t}_{ \underline r}(g)=0.$$

Set $\Gamma_{\underline r}:=\Gamma\cap  \underline  r$  equal to the pre-image of $\Lambda\cap \underline  r$ under the quotient $\Gamma\to \Lambda$. If $f$  is a function on
$\Gamma_{\underline r}$,  the hypotheses of Remark \ref{supconv} are satisfied and we can perform the convolutions  ${\mathcal P}^{\underline t\pm}_{ \underline r}*f$. So convolution by
${\mathcal Q}^{\underline t}_{ \underline r}$     induces a map
$$ \Pi^{\underline t}_{\underline r}:\mathcal C[\Gamma_{\underline r}]\to \mathcal C[\Gamma_{\underline t}],\quad f\mapsto {\mathcal Q}^{\underline t}_{ \underline r} *f.$$

Given a rational subspace $\underline r$,   let us consider the space $DM(X\cap \underline r)$ inside $\mathcal C[\Gamma_{\underline r}]$.

The following statement is similar to Proposition 3.6 in \cite{dpv}.
\begin{proposition}\label{indd}
$\Pi^{\underline t}_{\underline r}$ maps $DM(X\cap \underline r)$ to $DM(X\cap \underline t)$.
\end{proposition}
\begin{proof}

Let   $Y:=X\cap \underline r,\ Y':=X\cap \underline t$.
If $T\subset Y'$ is a cocircuit in $Y'$,
we need to see that $\nabla_T({\mathcal Q}^{\underline t}_{\underline r}* f)=0$
for each $f\in DM(X\cap \underline r)$.

By definition  $$\nabla_{A\cup B}{\mathcal P}^{\underline t+}_{ \underline r} =\nabla_{A\cup B}{\mathcal P}^{\underline t-}_{ \underline r} =
\delta_0$$ so that,
if $T=A\cup B=Y'\setminus Y$,
$$ \nabla_{T}(\mathcal Q^{\underline t}_{\underline r}*f)=
(\nabla_{T}{\mathcal Q}^{\underline t}_{\underline r})*f=0.$$
Otherwise, the set  $Y\cap T$ is a cocircuit in $Y$ and then $\nabla_{Y\cap T}f=0$. So
$$ \nabla_{T}({\mathcal Q}^{\underline t}_{\underline r}*f)=(\nabla_{T\setminus Y}{\mathcal Q}^{\underline t}_{\underline r})*(\nabla_{Y\cap T}f)=0.$$
\end{proof}

This proposition gives   a way to construct explicit elements of  $DM(X)$.
Take a flag  $\phi$ of oriented rational subspaces $0=\underline r_0\subset \underline r_1\subset  \underline r_2\subset \cdots  \subset \underline r_s=V$ with $\dim(\underline r_i)=i$.
Set
${\mathcal Q}_i:={\mathcal P}^{\underline r_{i}+}_{ \underline r_{i-1}}-{\mathcal P}^{\underline r_{i}-}_{ \underline r_{i-1}}$,
then
\begin{proposition}
\begin{equation}\label{lath} {\mathcal Q}_{ \phi}^X:={\mathcal Q}_1*{\mathcal Q}_2*\cdots*{\mathcal Q}_s \end{equation}
lies in $DM(X)$.
\end{proposition}

Set
$\theta_i^\pm(g):=\theta^{\underline r_{i}\pm}_{ \underline r_{i-1}}(g)$
the corresponding generalized functions on $G$, then their product is well defined.
We set
\begin{equation}\label{lath1} \theta_{ \phi}^X(g):=\prod_{i=1}^{s }( \theta^{ +}_i(g) -
\theta^{ -}_i(g)). \end{equation}

\subsection{Removing   a vector\label{rav}}

In this section, we prove the key technical result of this article, that is Theorem \ref{lasf}.

 Given a subgroup $\Psi\subset \Gamma$,  we identify  the space    $\mathcal C[\Gamma/\Psi]$ of $\mathbb Z$ valued functions on  $ \Gamma/\Psi$ with  the subspace of  $\mathcal C[\Gamma]$ formed by the functions constant on the cosets of $\Psi$. In particular, given $a\in  \Gamma,$  a function is constant on the cosets  of  $\mathbb Z a$ if and only if $\nabla_af=0$.  Therefore  the space  $\mathcal C[\Gamma /\mathbb Z a]$ is identified with $\ker(\nabla_a)$.

We define \begin{equation}\label{ilia}
i_a:\mathcal C[\Gamma /\mathbb Z a]\to \ker(\nabla_a),\ \quad i_a(f)(x):=  f(x+\mathbb Z a).
\end{equation}

The space
$DM(X\cap
\underline r )$ is a space of functions on  $\Gamma_{\underline r}:=\Gamma\cap \underline r$. We embed $\mathcal C[\Gamma_{\underline r}]$ into $\mathcal C[\Gamma]$ by extending each function    by  $0$ to the entire $\Gamma$ outside $\Gamma\cap \underline r$.

  \begin{definition}\label{DMG}
  Define $DM^{(G)}(X\cap
\underline r )$  to be  the $R(G)$ submodule of $\mathcal C[\Gamma]$  generated by the image of $DM(X\cap
\underline r )$ inside $\mathcal C[\Gamma]$.

\end{definition}
Let us assume now that $X$ is non degenerate.  Using Proposition \ref{indd}, we obtain,
 for each $\underline r\in S_X^{(s-1)}$,  a homomorphism
$$\Pi_{\underline r}:DM^{(G)}(X\cap \underline r)\to DM(X),\qquad \Pi_{\underline r}:=\Pi_{\underline r}^V$$
of $R(G)$ modules given by  convolution with  ${\mathcal Q}^{V}_{ \underline r}(X)$, the function given by Definition   \ref{defiQ}.

We can thus consider the map $$\Pi_X:H^X\to DM(X)$$ where
$$\Pi_X:=\oplus_{\underline r\in S_X^{(s-1)}}\Pi_{\underline r},\quad H^X:=\oplus_{\underline r\in S_X^{(s-1)}}DM^{(G)}(X\cap \underline r).$$
Take $a\in X$ of infinite order. Write $X=[Z,a]$ and set $\tilde Z$ to be the image of $Z$ in  $\Gamma /\mathbb Z a$.

\begin{theorem}\label{lasf}
1) We have an exact sequence:
\begin{equation}\label{labellissima}
0\to DM(\tilde Z )\stackrel {i_a}\to  DM(X )   \stackrel{\nabla_a}\to  DM( Z)\to 0.  \end{equation}

2) The map $$\Pi_X:H^X\to DM(X)$$
is surjective.
\end{theorem}
\begin{proof}
We proceed by induction on the number of elements in $X$.

 We know that the kernel of $\nabla_a$  in $\mathcal C [\Gamma]$  equals  $i_a(\mathcal C[\Gamma /\mathbb Z a])$. We first show that  $\ker(\nabla_a)\cap DM(X)=i_a(DM(\tilde Z))$.
This   is clear since a cocircuit in $\tilde Z$ is the image of a cocircuit in $X$ not containing $a$.

It is also  clear that $\nabla_a$ maps $DM(X)$ to $DM(Z)$. In fact, let   $T\subset Z$ be a cocircuit in $Z$,   then $T\cup   \{a\}$ is a cocircuit in $X$.  Thus $\nabla_T(\nabla_af)=0$ if $f\in DM(X)$.

Thus in order to finish the proof of the exactness of the sequence (\ref{labellissima}),  it only remains to show that the map
$\nabla_a:DM(X )    \to  DM( Z)$
is surjective.

If $\{a\}$ is a cocircuit, that is if  $Z$ is degenerate, then $DM(Z)=0$, the map $i_a$ is an isomorphism  and there is nothing more to be proven.

Assume that $Z$ is non degenerate. Given $\underline r\in S_X^{(s-1)}$, we have three possibilities.

In the first two cases,  we assume that $a\in \underline r$ so that  $X\cap \underline r=(Z\cap \underline r)\cup \{a\}$ and  set $\tilde {\underline r}$ equal to the hyperplane in $V/\langle a\rangle$ image of $\underline r$.
\begin{enumerate}
\item    $Z\cap \underline r$ does not span $\underline r$, that is    $\underline r\in S_X^{(s-1)}\setminus   S_Z^{(s-1)}$. In this case,   we get the isomorphism
$$DM^{(G)}(\tilde Z\cap \tilde {\underline r})\stackrel {i_a}\to DM^{(G)}(X\cap {\underline r}).$$
Furthermore $\Pi_{\underline r}\circ  i_{a}=i_a\circ \Pi_{\tilde {\underline r}}.$

\item   $Z\cap \underline r$   spans $\underline r$. In this case, by induction,  we may assume that we have the exact sequence
$$
0\to DM^{(G)}(\tilde Z\cap \tilde {\underline r} )\stackrel {i_{a}}\to  DM^{(G)}(X\cap {\underline r} )   \stackrel{\nabla_a}\to  DM^{(G)}( Z \cap {\underline r})\to 0.  $$
Furthermore $\Pi_{\underline r}\circ  i_{a}=i_a\circ \Pi_{\tilde {\underline r}}.$

\item  $a\notin \underline r$.  In this case, we get an equality $DM^{(G)}(X\cap {\underline r} ) = DM^{(G)}( Z \cap {\underline r})$.   Denote  by $\Pi_{\underline r}(X),\Pi_{\underline r}(Z)$ respectively the two maps  given  by convolution by the functions ${\mathcal Q}^{V}_{ \underline r}(X), {\mathcal Q}^{V}_{ \underline r}(Z)$ associated to the two lists $X$ and $Z$.  We clearly have $\nabla_a{\mathcal Q}^{V}_{ \underline r}(X)= {\mathcal Q}^{V}_{ \underline r}(Z)$ and hence  the identity
$\Pi_{\underline r}(Z)=\nabla_a\Pi_{\underline r}(X)$.
\end{enumerate}

Define a map  $p:H^X \to  H^Z $ as follows.  If   $\underline r\in S_X^{(s-1)}\setminus   S_Z^{(s-1)}$,  we are in  case  i) and we set $p=0$ on $DM^{(G)}(X\cap {\underline r} )$.    In case ii)   $\underline r\in     S_Z^{(s-1)}$  and $a\in \underline r$, we set $p:=\nabla_a:DM^{(G)}(X\cap {\underline r} ) \to DM^{(G)}( Z \cap {\underline r})$,
a surjective map by induction. Finally in case iii)   $\underline r\in     S_Z^{(s-1)}$  and $a\notin \underline r$,  then  $DM^{(G)}(X\cap {\underline r} )= DM^{(G)}( Z \cap {\underline r})$  and we take $p=Id$ the identity.

We have thus that $p$ is surjective and its kernel is the  direct sum of $DM^{(G)}(X\cap {\underline r} )$ for  $\underline r\in S_X^{(s-1)}\setminus   S_Z^{(s-1)}$ plus  the direct sum
of the kernels of the maps
$\nabla_a:DM^{(G)}(X\cap {\underline r} ) \to DM^{(G)}( Z \cap {\underline r})$,  over the $\underline r\in     S_Z^{(s-1)}$ with  $a\in \underline r$.

Let us next compare   $\ker (p)$ with $ H^{\tilde Z}$.  The preimage in $V$  of a  subspace  $\underline v\in S_{\tilde Z}^{(s-2)}$ is  a subspace $\underline r\in  S_X^{(s-1)}$ with $a\in\underline r$,  thus  a space of the first two types. By all the previous remarks,  taking direct sums and assuming by induction on dimension our result true for each list $X\cap\underline r$, $\underline r\in S_X^{(s-1)}$, we deduce an exact sequence
$$ 0\to H^{\tilde Z}\to H^X\to H^Z\to 0.$$ Furthermore the diagram
\begin{equation}\label{ildiafr}\begin{CD}\hskip-0.7cm 0\to H^{\tilde Z}@>>>H^X@>>>\hskip0.7cm  H^Z\to 0\\  @V{\Pi_{\tilde Z}}VV@V{\Pi_X}VV@V{\Pi_Z}VV@.\\
0\to DM(\tilde Z )@>{i_a}>>  DM(X )  @>{\nabla_a}>>  DM( Z)\to 0\end{CD}\end{equation}
commutes.

By induction,  the map  $\Pi_Z$ is surjective. We immediately deduce the surjectivity of $\nabla_a:DM(X)\to DM(Z)$,  hence part 1).

2) To finish,  we have to prove the surjectivity of $\Pi_X$. By induction,  both $\Pi_{\tilde Z}$ and $\Pi_{Z}$ are surjective so that also $\Pi_X$ is,  as desired.
\end{proof}

Let $s:=\dim V$, we say that a sublist $\underline b:=\{b_1,\ldots,b_s\}$ of $X$ is {\it a basis}  if the elements  $\overline{\underline b}:=\{\overline b_1,\ldots,\overline b_s\}$ form a basis for $V$.  Let us denote by $\mathcal B_X$ the set of all such bases extracted from $X$.

For each    $\underline b\in \mathcal B_X$,    the elements  $\{\overline b_1,\ldots,\overline b_s\}$ generate a lattice  $\langle \overline b_1,\ldots,\overline b_s \rangle$ in $\Lambda$ of some index $d(\underline b)$  (equal to the absolute value of the determinant or the volume of the corresponding parallelepiped).

The elements $\underline b$ generate a free abelian group $\mathbb Z_{\underline b}$  and  the subgroup  $\Gamma_t\times \mathbb Z_{\underline b}$  is still of index $d(\underline b)$ in $\Gamma$.
Set
$$\delta(X):=\sum_{\underline b\in \mathcal B_X}d(\underline b).$$
The number $\delta(X)$ has a geometric meaning,  it is the volume of the {\it zonotope}  $B(X):=\{\sum_{a\in X}t_a\overline a\,|\, 0\leq t_a\leq 1\}$.  This   is  a consequence of the decomposition of the zonotope into parallelepipeds  of volume $d(\underline b)$ (cf. \cite{dp1}).
\begin{corollary}
$DM(X)$  is a free $\mathbb  Z [\Gamma_t]$-module of rank $\delta(X)$.
\end{corollary}
\begin{proof}
Notice that, if $\Gamma$ is finite or $X$ consists of torsion elements,  our statement is trivially true.

We  now proceed by induction  using the exact sequence (\ref{labellissima}).  Thus, take $a\in X$ an element of infinite order and write $X=[Z,a]$. By induction,  the space   $ DM( Z)$  is a free  $\mathbb  Z [\Gamma_t]$-module of rank $\delta(Z)$,  thus it suffices to see that $DM(\tilde Z )$ is a free $\mathbb  Z [\Gamma_t]$-module of rank $\delta(X)-\delta(Z).$

Set $\tilde \Gamma= \Gamma/\mathbb  Z a$. By induction,   $DM(\tilde Z )$ is a free $\mathbb  Z [\tilde \Gamma_t]$ module of rank $ \delta(\tilde Z).$  On the other hand, since $a$ is not a torsion element,   $\mathbb  Z [\tilde \Gamma_t]$ is a free   $\mathbb  Z [\Gamma_t]$-module of rank  $|\tilde \Gamma_t/ \Gamma_t|$. Thus it suffices to see that
$\delta(X)-\delta(Z)=|\tilde \Gamma_t/ \Gamma_t| \delta(\tilde Z).$

  The number $\delta(X)-\delta(Z) $ is the sum, over all the bases $\underline b$ of  $X$ containing $a$, of  the index   $d(\underline b)$ of $\Gamma_t\times \mathbb Z_{\underline b}$    in $\Gamma$. Pass  modulo  $a\in \underline b$  and denote by $\tilde { \underline b}$  the induced basis in $\tilde \Gamma $. We see  that $$|\Gamma_t|d( { {\underline b}})=| \Gamma/ \mathbb Z_{\underline b}|=| \tilde \Gamma/\mathbb Z_{\tilde {\underline b}}|=|\tilde \Gamma_t | d(\tilde {\underline b})$$
as desired.
\end{proof}
\begin{remark}
One could have also followed the approach of Dahmen--Micchelli to this theorem. Then one needs to show that  the space $DM(X)$  restricts isomorphically to the space of $\mathbb Z$ valued functions on any set $\delta(u\,|\,X):=\Gamma\cap ( u-B(X))$  when  $u$ is a generic
 vector. In turn,  this set is a union of $\delta(X)$ cosets of $\Gamma_t$  and this statement can be interpreted as saying that:

{\it The values on $\delta(u\,|\,X)$  are initial values for the recursion equations given by the cocircuits.  A function in $DM(X)$ is completely determined by these values that can be assigned arbitrarily } (see \cite{dp1}).
\end{remark}
Proceeding by induction and using the surjectivity of $\Pi_X:H^X\to DM(X)$ we immediately get

\begin{corollary}\label{geni}
The elements $\mathcal Q_{ \phi}^X$, as $\phi$ varies among complete flags of rational subspaces, span $DM(X)$ as $R(G)$ module.
\end{corollary}

\subsection{A basis for $DM(X)$}  This section is not needed for the paper but it gives a very precise description of $DM(X)$.\smallskip

We assume that $X$ is totally ordered in such a way
that the torsion elements follow the elements of infinite order.

Take a basis $\underline b\in \mathcal B_X$.  Given $c\in X$, set $\underline b_{\geq c}=\{b\in \underline b\,|\,c\leq b\}$. Define $X_{\underline b}=\{c\in X\,|\,\overline c\in\langle \underline b_{\geq c}\rangle\}$.
Since $X_{\underline b}$ is a sublist of $X$,  we have that  $DM(X_{\underline b})$ is a subset of  $DM(X)$.

Now consider the complete flag $\phi_{\underline b}$ of rational subspaces $\underline r_i$ generated by $ \{\overline b_1,\ldots,\overline b_i\},\ i=1,\ldots,s$.
We obtain the element  ${\mathcal Q}^{X_{\underline b}}_{\phi_{\underline b}}\in DM(X_{\underline b})\subset DM(X)$.

Choose a complete set $C_{\underline b}\subset \Gamma$ of coset representatives modulo the subgroup  $\Gamma_t\times \mathbb Z_{\underline b}$.
To every pair $(\underline b,\gamma)$ with $\underline b\in \mathcal B_X$ and $\gamma\in C_{\underline b}$, we associate the element $${\mathcal Q}^{X_{\underline b}}_{\phi_{\underline b},\gamma}=\tau_\gamma {\mathcal Q}^{X_{\underline b}}_{\phi_{\underline b}}$$
in  $DM(X).$ We can now state  our

\begin{theorem}\label{LABAS} The set of elements ${\mathcal Q}^{X_{\underline b}}_{\phi_{\underline b},\gamma}$, as $\underline b$ varies in $\mathcal B_X$ and $\gamma$ in $C_{ \underline b } $, is a basis of $DM(X)$ as a $\mathbb Z[\Gamma_t]$ module. \end{theorem}
\begin{proof}
Let $a$ in $X$ be the least element. As usual,  we can assume that it is of infinite order. As in the previous section,  write $X=[a,Z]$ and set $\tilde Z$ to be the image of $Z$ in  $\tilde \Gamma:=\Gamma /\mathbb Z a$.

Take a basis $\underline b\in \mathcal B_X$. If $a$ does not lie in $\underline b$, then $\underline b\in \mathcal B_Z$. Observe that $Z_{\underline b}=X_{\underline b}\setminus \{a\}$ and that for every $\gamma\in C_{\underline b}$
$${\mathcal Q}^{Z_{\underline b}}_{\phi_{\underline b},\gamma}=\nabla_a({\mathcal Q}^{X_{\underline b}}_{\phi_{\underline b},\gamma}).$$

Thus, by induction, the elements $\nabla_a({\mathcal Q}^{X_{\underline b}}_{\phi_{\underline b},\gamma})$, as $\underline b$ varies in $\mathcal B_Z$ and $\gamma$ in $C_{ \underline b } $,  are a basis of $DM(Z)$ as a $\mathbb Z[\Gamma_t]$ module.

If on the other hand $a$ lies in $\underline b$, consider the corresponding basis $\tilde {\underline b}\in  \mathcal B_{\tilde Z}$.
Take  the quotient map $p:\Gamma\to \tilde \Gamma$. It is clear that $p$ maps $C_{\underline b}$ injectively onto a set of representatives of the cosets of $\Gamma_t\times  \mathbb Z_{\tilde{\underline b}}$ in $\tilde \Gamma$.
Define  an equivalence relation on $p(C_{\underline b})$ setting $\gamma\simeq \gamma'$ if $\gamma - \gamma' \in \tilde\Gamma_t\times  \mathbb Z_{\tilde{\underline b}}.$ If we choose an element in each equivalence class,  we obtain a subset of $\tilde\Gamma$ which we can take as $C_{\tilde{\underline b}} $.

By induction,  it is then clear that the set of elements ${\mathcal Q}^{\tilde Z_{\tilde{\underline b}}}_{ \phi_{\tilde {\underline b}}, p(\gamma )}$, as $\tilde{\underline b}$  varies in $\mathcal B_{\tilde Z}$  and  $\gamma$ in $C_{ { \underline b} } ,$ is  a basis of $DM(\tilde Z)$ as a $\mathbb Z[\Gamma_t]$ module.
Now observe    that

\begin{equation}\label{formulac}
i_{a}({\mathcal Q}^{\tilde Z_{\tilde{\underline b}}}_{\phi_{\tilde{ \underline b}},p(\gamma)})={\mathcal Q}^{X_{\underline b}}_{\phi_{\underline b},\gamma}.
\end{equation}
 In fact, whenever one multiplies a function  $f=i_a(g)\in \ker \nabla_a$
 by any element $ \delta_c $ (or by a series formed by these elements), one sees that
\begin{equation}\label{aal}
 \delta_c * i_a(g)=i_a ( \delta_{p(c)}*  g),\ \forall c\in\Gamma,\ \forall g\in \mathcal C[\tilde \Gamma].
\end{equation}
Furthermore the elements  $(X_{  \underline b } \setminus\{a\})\cap \underline r_i,\ i\geq 2,$ map surjectively  under $p$  to the elements  of  $\tilde Z_{p(\underline b)}\cap  \underline r_i/ \underline r_1$  and thus the formulas defining the two functions on both sides of Equation (\ref{formulac})
coincide, factor by factor, for all steps of the flag  from Formula (\ref{aal}).

Thus everything follows from the exact sequence (\ref{labellissima}).
\end{proof}

\subsection {Support of $DM(X)$. }\label{supports}
In this subsection, we are going to assume that $X$ spans $V$. Otherwise $DM(X)=\{0\}$ and our discussion is trivial.

Recall that we have denoted by $\mathcal B(X)$ the set of all bases extracted from $X$.
Let $\underline b=[b_1,b_2,\ldots, b_s]$ be in  $\mathcal B_X$. Correspondingly we get an inclusion
$j_{\underline b}:\mathbb Z^s\to \Gamma$ defined  by $j_{\underline b}((n_1,\ldots ,n_s))=\sum_{h=1}^sn_hb_h$ and  a surjective homomorphism $j^*_{\underline b}:G\to (S^1)^s$. We denote its kernel by $G(\underline b)$. This is a finite subgroup in $G$.
We define

 \begin{definition}
 \begin{enumerate}
 \item
  $P(X)$ is  the union of all the sets $G(\underline b)$, when
 $\underline b$ varies in  $\mathcal B_X$.
 \item
For $g\in G$, denote by $X_g$ the sublist of elements of $X$ taking value $1$ at  $g$.
\end{enumerate}
\end{definition}

Now consider the ring $R(G)$ as a ring of functions on $G$, or better its complexification $G_{\mathbb C}$.

\begin{lemma}\label{lusuppo}  Let $V(X)$  denote the set of
  elements $g\in G_{\mathbb C}$  such that  $\prod_{a\in Y}(1-g^a)=0$,
when $Y$ runs over all cocircuits. Then $V(X)= P(X)$.
\end{lemma}
\begin{proof}
If $g\in G(\underline b)$ and $Y$ is a cocircuit, necessarily there exists an element $b_i$ of $\underline b$ belonging to $Y$. Otherwise the complement of $Y$ would not be contained in a proper subspace of $V$. So we have $\prod_{a\in Y}(1-g^a)=0$.

 Conversely take $g\in V(X)$, and consider $X_g$. We claim that $X_g$     spans $V$, otherwise $X\setminus X_g$ contains a cocircuit $Y$,  and $\prod_{a\in Y}(1-g^a)\neq 0$. So there is a basis $\underline b$ in $X_g\subset X$ and therefore  $g\in G(\underline b)$
 \end{proof}

Since the ideal $J_X$ defines the finite set of points   $P(X)$, we have  that $\mathbb C[\Gamma]/J_X $ is a semi--local algebra,  direct sum of its local components at each point in $  P(X)$.

  Let us normalize the Haar measure on $G$ to be of total mass $1$. This allows us to identify generalized functions on $G$ and distributions on $G$. Under this identification,  call $\widehat {DM} (X)$ the space of distributions  on $G$ of which $DM(X) $ gives the Fourier coefficients. Thus,
if $p\in G$, the  delta distribution at $p$ is identified to the generalized function $$\delta_{p}(g)=\sum_{\gamma\in \Gamma}p^{-\gamma}g^{\gamma},$$
the Fourier transform of  the element $f_p\in \mathcal C[\Gamma]$ given by
 $f_{p}(\gamma):=p^{-\gamma}.$

Any generalized function on $G$ supported at $p$ is a derivative of the $\delta$ function $\delta_{p}$, so is
of the form
$$\theta(p,q)= \sum_{\gamma\in \Gamma}p^{-\gamma} q( \gamma)g^\gamma,$$
 where $q$ is a polynomial on $V$. We have written $q(\gamma)$ instead of $q(\overline \gamma)$, where $\overline \gamma$ is the image of $\gamma$ in $V$.
The function $\theta(p,q)$  is
 the Fourier transform of
$$f(p,q)(\gamma)= p^{-\gamma} q(\gamma) .$$

 From the definitions,   $DM_{\mathbb C}(X)$   is dual to the   algebra $\mathbb C[\Gamma]/J_X $
  thus, from Lemma \ref{lusuppo},  we deduce
 \begin{proposition}\label{lesuppo} The support of every element in $\widehat {DM} (X)$  is contained in the finite set $P(X)$.
 \end{proposition}Moreover,  once we take complex linear combinations of the elements in $\widehat {DM} (X)$,  we obtain a space of distributions which is the direct sum  of
  the local contributions to $\widehat {DM}(X)$ at each point $p\in P(X)$. In order to understand these local contributions,   we need to
recall \cite{dpv} the  {\it differentiable} Dahmen-Micchelli space $D(X)$.

 Given a vector $v\neq 0$ in $V$,
 we denote by $\partial_v$ the directional
derivative associated to $v$. For a subset $Y$ of $X$ ,
we denote   by $\partial_{Y}:=\prod_{a\in Y}\partial_{\overline a}$.

\begin{definition}
The space $D(X)$ is the space of polynomial functions on $V$
   satisfying the system of
differential equations $\partial_Yf=0$ as $Y$ varies among all cocircuits
of $X$.
\end{definition}

\begin{remark} The space $D(X)$ is   a space of polynomials of finite dimension $d(X)$:  the number of bases of $V$ extracted from $X$. It describes locally some important functions as the box--spline.

We will prove in a subsequent article that $D(X)$ is isomorphic to the $G$-equivariant cohomology   with compact support of $M_X^f$. Similarly, for each $g\in P(X)$,   the space $D(X_g)$  is the $G$-equivariant cohomology  with compact support of the submanifold  $(M_X^f)^g$  consisting  of  fixed points by $g$ in $M_X^f$.

The equivariant cohomology, as an algebra, is   the algebra of differential operators with constant coefficients induced on $D(X)$.
\end{remark}

The following proposition is proven in \cite{DM}.
\begin{proposition}
Any function $f$ in $DM_{\mathbb C}(X)$ can be written uniquely as
$$f=\sum_{p\in P(X)} f(p,q)$$
where $q\in D(X_p)$.
\end{proposition}
\begin{proof}
We recall briefly the proof.
We have just seen  that $f\in DM(X)$ can be written uniquely as $f=\sum f(p,q)$
, where $p\in P(X)$ and $q$ is a polynomial on $V$.

We verify that $\nabla_Y f(p,q)=f(p, \nabla(Y,p)q)$, where $\nabla(Y,p)=\prod_{a\in Y}(1-p^a \tau_a)$.
On polynomials, the  operator $(1-\tau_a)q(v)=q(v)-q(v-a)$ is nilpotent, so that
$(1-p^a \tau_a)=(1-p^a)-(p^a(\tau_a-1))$ is invertible when $p^a\neq 1$.

Thus we split $X$ in $X_p\cup (X\setminus X_p)$ and  see that  $\nabla_Y f(p,q)=0$ if and only $\nabla_{Y\cap X_p}q=0$.
If $Y$ is a cocircuit of $X$, then $Y\cap X_p$ contains a cocircuit of $X_p$. Thus
$f(p,q)\in DM(X)$, if and only if $\nabla_{ Z} q=0$ for all cocircuits $Z$ of $X_p$.
Similarly, if $v\in V$,  the difference operator $\nabla_v$  and the differential
operator  $\partial_v$
satisfy  $\nabla_v=T_v\partial_v$ with $T_v=\frac{1-e^{-\partial_v}}{\partial_v}$   invertible on the space of polynomials.
Thus  we obtain that $f(p,q)\in DM_{\mathbb C}(X)$  if and only if $q\in D(X_p)$.
\end{proof}
\begin{remark}
Since each point $p\in P(X)$ is an element of finite order in $G$,  we clearly see the quasi--polynomial nature of each summand  $f(p,q)$.
\end{remark}

\section{The spaces $\mathcal F(X)$ and $\tilde{\mathcal F}(X)$}

Following  \cite{dpv}, we define

\begin{definition}
The space $\mathcal F(X)$  is the space of functions  $f\in \mathcal C[\Gamma]$  such that $\nabla_{X\setminus  \underline r }
f \text{ is supported on }  \underline r \text{ for every
proper rational subspace    } \underline r.
$
\end{definition}

Notice that, since $DM(X)$ is  the space of integer valued functions on $\Gamma$ such that $\nabla_{X\setminus  \underline r }
f=0$ for each rational subspace $\underline r$, $DM(X)\subset \mathcal F(X)$.

In \cite{dpv}, we have introduced  the space  $\mathcal F(X)$  only when $X$ is a finite list of characters of a connected torus $G$,  but the discussion there repeats verbatim in this general case.
We have a {\it canonical filtration}  which we shall interpret geometrically in the next subsection.
Consider   $\nabla_{X\setminus\underline r}$ as an operator on
$\mathcal F(X)$ with values in ${\mathcal C}[\Gamma]$.
Define the spaces
$$\mathcal F_{i}(X):=\cap_{\underline t\in S_X^{(i-1)}}\ker
\nabla_{X\setminus \underline t}\cap \mathcal F(X).$$

Notice that   by definition: $\mathcal F_{0}(X)=\mathcal F(X)$,
$\mathcal F_{\dim V}(X) $ is the space  $DM(X)$ and   $ \mathcal F_{i+1}(X) \subset
\mathcal F_{i}(X)$. It follows from the definitions that, for each $\underline t\in S_X^{(i)}$, $\nabla_{X\setminus \underline t}$ maps $\mathcal F_{i}(X)$ to $ DM(X\cap
\underline t )$  and $\mathcal F_{i+1}(X)$ to $0$. Consider the map
$$\mu_i=(\oplus\nabla_{X\setminus \underline t})_{\underline t\in S_X^{(i)}}:\mathcal F_{i}(X)\to \bigoplus_{\underline t\in S_X^{(i)}}
DM(X\cap
\underline t ).$$

With a proof entirely similar to the proof of Lemma 3.10  in \cite{dpv},  we obtain the following theorem.
\begin{proposition}\label{split} The sequence
$$0\to \mathcal F_{i+1}(X)\to\mathcal F_{i}(X)\stackrel{\mu_i}\to\bigoplus_{\underline r\in S_X^{(i)}}
DM(X\cap
\underline r )\to 0$$
is split exact with   splitting:
$$\mathcal F_{i}(X)=\mathcal F_{i+1}(X)\oplus\bigoplus_{\underline r\in S_X^{(i)}}
{{\mathcal P}}_{X\setminus \underline r}^{F_{\underline r} }*DM(X\cap
\underline r ).$$
\end{proposition}
As a consequence we obtain a decomposition dependent  upon the choices of the faces $F_{\underline r}$\begin{equation}\label{decF}
\mathcal F(X)=DM(X)\oplus \Big(\oplus_{\underline r\in S_X| {\underline r} \neq V}
{{\mathcal P}}_{X\setminus \underline r}^{F_{\underline r}}*DM(X\cap
\underline r )\Big).
\end{equation}
\begin{remark}
From Formula (\ref{decF}) and Theorem \ref{LABAS},   one can describe an explicit linear basis for the abelian group  $\mathcal F(X)$.
\end{remark}

While $DM(X)$ is a $R(G)$ module,
in general   $\mathcal F(X)$ is not stable under $R(G)$,
so we  consider the  $R(G)$  submodule $\tilde{\mathcal F}(X)$  in $\mathcal C[\Gamma ]$ generated by $\mathcal F(X)$ and similarly for   $\tilde{\mathcal F}_i(X)$.
\begin{remark}
It is easy to see that, if $X'$ is deduced from $X$ as in Remark \ref{real},
  then  $\tilde{\mathcal F}(X)=\tilde{\mathcal F}(X')$.
\end{remark}

 Recall the definition \ref{DMG} of the subspace $DM^{(G)}(X\cap
\underline r )$  generated by $DM(X\cap
\underline r )$ inside $\mathcal C[\Gamma]$.

      The group  $\Gamma_{\underline r}=\Gamma\cap \underline r$ is a direct summand in $\Gamma$ containing $\Gamma_t$.
       Thus the group $G_{\underline r}$, kernel of all the elements in $\Gamma_{\underline r}$, is connected  with character group $\Gamma/ \Gamma_{\underline r}$.
One gets a split exact sequence
$$1\to G_{\underline r}\to G\to G/G_{\underline r}\to 1.$$
We deduce that
$$   R(G/G_{\underline r} )=\mathbb Z[\Gamma_{\underline r}]$$
while $ R(G)$ is isomorphic (in a non canonical way) to
$$\mathbb Z[\Gamma]\sim \mathbb Z[\Gamma_{\underline r}]\otimes \mathbb Z[\Gamma/\Gamma_{\underline r}]=R(G/G_{\underline r})\otimes R(G_{\underline r}) .$$

It is then immediate to verify that \begin{equation}\label{IDM}
DM^{(G)}(X\cap
\underline r )=R(G)\otimes_{R(G/G_{\underline r})}DM(X\cap
\underline r )=R(G_{\underline r})\otimes DM(X\cap
\underline r ).
\end{equation}

Since both the maps $\mu_i$ and the convolutions ${\Theta}_{X\setminus \underline r}^{F_{\underline r}}$ obviously extend to maps of $R(G)$ modules, we easily deduce
\begin{corollary}\label{spit2} The sequence
$$0\to \tilde {\mathcal F}_{i+1}(X)\to\tilde{\mathcal  F}_{i}(X)\stackrel{\mu_i}\to\bigoplus_{\underline r\in S_X^{(i)}}
DM^{(G)}(X\cap
\underline r )\to 0$$
is split exact with   splitting:
$$\tilde {\mathcal F}_{i}(X)=\tilde {\mathcal F}_{i+1}(X)\oplus\bigoplus_{\underline r\in S_X^{(i)}}
{\mathcal P}_{X\setminus \underline r}^{F_{\underline r}}*DM^{(G)}(X\cap
\underline r ).$$
\end{corollary}
From this,  we obtain a decomposition dependent  upon the choices of the faces $F_{\underline r}$\begin{equation}\label{decFt}
\tilde{\mathcal F}(X)=DM(X)\oplus \Big(\oplus_{\underline r\in S_X| {\underline r} \neq V}
{{\mathcal P}}_{X\setminus \underline r}^{F_{\underline r}}*DM^{(G)}(X\cap
\underline r )\Big).
\end{equation}
 One easily   verifies  that $\tilde{\mathcal F}(X)$ can also be intrinsically defined as:
\begin{definition}
The space $\tilde{\mathcal F}(X)$  is the space of functions  $f\in \mathcal C[\Gamma]$  such that $\nabla_{X\setminus  \underline r }
f$  is supported on a finite number of $\Gamma$ translates of  $ \underline r$  for every
proper rational subspace    $ \underline r.
$
\end{definition}

\begin{remark}
By definition,  the filtration  $\tilde{\mathcal F}_i(X)$  is defined by a condition on torsion   or in the language of modules by support.  One can verify that in particular  $DM(X)$  is the part supported  in dimension 0, that is of maximal possible torsion.
\end{remark}

\subsection{Generators of $\tilde{\mathcal F}(X)$\label{genf}}

In this subsection, we assume that  $\overline a\neq 0$  for any $a\in X$.
Thus every open face $F$ produces a decomposition $X=A\cup B$ into positive and negative  vectors  and can define as  in  (\ref{esplicit}):
$$\mathcal P_X^{F}:=(-1)^{|B|}\delta_{-b_B}* H_A*H_{-B}.$$

\begin{theorem}\label{GENF}
The elements  $\mathcal P_X^{F}$, as $F$ runs on all open faces, generate  $\tilde{\mathcal F}(X)$ as $R(G)$ module.
\end{theorem}
\begin{proof}
 Denote by $M$ the $R(G)$ module generated by the elements $\mathcal P_X^{F}$, as $F$ runs on all open faces.   In general,  from  the description of $\tilde{\mathcal F}(X)$ given in
 \ref{split},  it is enough to prove that  elements of the  type   $
{\mathcal P  }_X^{F_{\underline r}}*g$   with $g\in DM(X\cap
\underline r)$ are in $M$. As $DM(X\cap \underline r)\subset \mathcal F(X\cap \underline r)$,
it is sufficient to prove by induction that each element
${\mathcal P}_X^{F_{\underline r}}*{\mathcal P}_{X\cap \underline r}^{K}$
is in $M$,  where $K$  is any open face for the system  $X\cap
\underline r$.
We choose a   linear function $u_0$ in the face $F_{\underline r}$. Thus $u_0$   vanishes on  $\underline r$ and  is non zero on every element $a\in X$ not in $\underline r$. We choose a linear function $u_1$ such that the restriction of $u_1$ to $\underline r$ lies in the face $K$.  In particular, $u_1$ is non zero on every element $a\in X\cap \underline r$.
 We can choose $\epsilon$ sufficiently small such that $u_0+\epsilon u_1$ is  non zero on every element $a\in X$. Then $u_0+\epsilon u_1$ defines an open face $F$ in the arrangement $\mathcal H_X$.  We see that ${\mathcal P  }_X^{F_{\underline r}}*{\mathcal P}_{X\cap \underline r}^{K}$ is equal to ${\mathcal P}_X^F$.
\end{proof}\bigskip

\subsection {Support of $\tilde{\mathcal F}_i(X)$.}    For a fixed $0\leq i\leq s$,
 take a   linearly independent sublist
$\underline [b_1,\ldots ,b_i]$ in $X$.
Correspondingly we get an inclusion
$j_{\underline b}:\mathbb Z^i\to \Gamma$ defined  by $j_{\underline b}((n_1,\ldots ,n_i))=\sum_{h=1}^hn_hb_h$ and  a surjective homomorphism $j^*_{\underline b}:G\to (S^1)^i$. We denote its kernel by $G(\underline b)$. This is a   subgroup in $G$ of dimension $s-i$.
We define

 \begin{definition}
  $P_{s-i}(X)$ is  the union of all the sets $G(\underline b)$, when
 $\underline b$ varies among linearly independent sublists in $X$ of cardinality $i$.
\end{definition}

As before,  normalize the Haar measure on $G$ to be of total mass $1$ and  identify generalized functions on $G$ and distributions on $G$.  Call
 $\widehat {\tilde{\mathcal F}}_i (X)$ the space of distributions  on $G$ of which $\tilde{\mathcal F}_i(X)$.   gives the Fourier coefficients.

Then the decomposition (\ref{decFt}) together with Proposition \ref{lesuppo} immediately implies

 \begin{proposition}\label{supportFi} The support of every element in $\widehat {\tilde{\mathcal F}}_i (X)$  is contained in   $P_{s-i}(X)$.
 \end{proposition}

\section{Index theory\label{inteo}}

\subsection{$K$-theory}

We briefly  review the  notations for  $K$-theory that we will use.

Let $G$ be a compact Lie group acting  on a locally compact  space $N$.
One  has the notion of the equivariant  topological $K-$ theory group $K_G^0(N)$.  $K_G^0(N)$ is a contravariant functor for proper maps and covariant for   open embeddings.
We recall that  representatives of the $K-$theory group $K_G^0(N)$  can  be described in the following way. Given two  $G$-equivariant complex vector bundles  $E^0,E^1$ on $N$
  and a $G$-equivariant bundle  map  $f:E^0\to E^1$, the support $supp(f)$ of $f$ is the set of points where $f_x:E^0_x\to E^1_x$ is not an isomorphism.
  A $G$-equivariant bundle map $f$ with compact support   defines an element $[f]$ of $K_G^0(N)$.

 Let $f:E^0\to E^1$ and $g: F^0\to F^1$ be two $G$-equivariant bundle maps.
 Using $G$-invariant  Hermitian metrics on the bundles $E^i,F^i$ we can define:  $$f\odot g:E^0\otimes F^0\oplus E^1\otimes F^1\to E^1\otimes F^0\oplus E^0\otimes F^1$$
by $$f\odot g:=\left(
                \begin{array}{cc}
                  f\otimes 1 & -1\otimes g^* \\
                  1\otimes g& f^*\otimes 1 \\
                \end{array}
              \right).$$

 The support of $f\odot g$ is the intersection of the supports of $f$, $g$  thus  $f\odot g$ induces  an element in $K_G^0(N)$ as soon as  one of the two $f,g$ has compact support.

 In particular this defines a product $[f][g]:=[f\odot g]$  on $K_G^0(N)$.

If $N=\rm pt$ is a point, then $K_G^0({\rm pt})$ is isomorphic  to the  Grothendieck ring  $R(G)$  of finite dimensional representations of $G$.

In general,  tensor product with finite dimensional  representations of $G$ induces a $R(G)$ module structure on $K_G^0(N)$. Take the projection $\pi:N\to \rm pt$.  Given $\tau\in R(G)$ and  $\sigma\in K_G^0(N)$,
  we have that $[\pi^*(\tau)\odot \sigma]\in K_G^0(N)$ and this gives a $R(G)$ module structure to $K_G^0(N)$.

We will need also the groups  defined inductively as $$K_G^{i+1}(N):=K_G^{i}(N\times\mathbb  R).$$ One has $K_G^0(\mathbb R)= K^1_G(\rm pt)=0$.

There is a natural isomorphism
$K_G^{i}(N)\to K_G^{i+2}(N)=K_G^i(N\times \mathbb R^2)$   given by Bott  periodicity, that we describe  below.

$\bullet$ Let $W$ be a Hermitian vector space and let $E=\bigwedge W$.
Then, for $w\in W$, consider the exterior multiplication  $m(w):E\to E$  and the Clifford  action
\begin{equation} \label{cli}
c(w)=m(w)-m(w)^*,\qquad (m(w)(\omega):=w\wedge \omega)
\end{equation}
of $W$ on $\bigwedge W$. Then one has
$c(w)^2=-\|w\|^2,$
so that $c(w)$ is an isomorphism, if $w\neq 0$.

\bigskip

If $p:W\to M$ is a $G-$equivariant complex vector bundle over a $G-$space $M$,
we can consider also $W$ as $G-$space and we have  a Thom  isomorphism
   $${C}_W: K_G^0(M)\to K_G^0( W).$$
  In order to make this explicit,  we use a $G$-invariant Hermitian metric on $W$.
  Then the fiberwise Clifford action   $c(w_x):\bigwedge^{even} W_x\to \bigwedge^{odd} W_x$ defines
  a morphism ${\bf c}_W: p^* \bigwedge^{even} W\to p^*\bigwedge^{odd}W$ of  vector bundles
  over $   W$, that we call the  {\it Bott symbol}.
   Take    a bundle map  $f:E\to F$ of complex equivariant vector bundles   on $M$ which is an isomorphism outside a compact set, and denote still by $f$  its pull back $f:p^*E\to p^*F$.
   Then  $f\odot {\bf c}_W$ is a bundle map of bundles over $  W$,  which is an isomorphism outside the support of $f$ embedded in $ W$ via the zero section. We set
   \begin{equation}
    { C}_W([f])=[f\odot {\bf c}_W].
    \end{equation}

$\bullet$
If $F$ is a $G$-invariant closed subset of $N$,  denote by $i:F\to N$  the closed embedding and $j:N\setminus F\to N$ the open  embedding. There is a long exact sequence
of $R(G)$ modules:
\begin{equation}\label{long}\cdots
\to K_G ^i(N\setminus F)\stackrel{j_*}\to K_G^{i}(N)\stackrel{i^*}\to K_G^i(F)\stackrel{\delta}\to K_G^{i+1}(N\setminus F)\to\cdots
\end{equation}

\subsection{Transversally elliptic operators}

To define the index of a transversally elliptic  operator on $M$, we require the following hypothesis on $M$.
If $G$ is a compact Lie group acting on $M$, we assume that $M$  can be embedded  as a $G$-invariant open subset of  a compact $G$-manifold $\tilde M$.

Given  such a  manifold $M$  with a $C^\infty$ action of a compact Lie group  $G$,
one has the notion of {\it transversally elliptic operator}  between two equivariant complex vector bundles  $E ,F$ on $M$. Such an operator is a pseudo-differential operator  $A:\Gamma(M,E)\to \Gamma(M,F)$ from  the space $\Gamma(M,E)$ of  smooth sections of $E$ to  the space $\Gamma(M,F)$ of  smooth sections of $F$, which commutes with the action of $G$,  is {\it elliptic} in the directions transversal to the orbits  of $G$ and is ``trivial '' at infinity.

 In more technical terms, let $T^*M$ denote the cotangent bundle of $M$ and $p:T^*M\to M$ the canonical projection. Inside $T^*M$, there is a special closed subset  denoted by $T^*_GM$.  Its fiber over a point $x\in M$ is  formed by all the cotangent vectors $\xi\in T^*_xM$  which vanish on the tangent space to the orbit of $x$  under $G$, in the point $x$. Thus   each fiber $(T^*_GM)_x$ is a linear subspace  of $T_x^* M$. In general the dimension of $(T^*_GM)_x$  is not constant and this space is not a vector bundle.
\begin{definition}
By a {\it symbol},  one means   a  smooth section  on $T^*M$  of  the bundle $\hom(p^*(E ),p^*(F))$:  in other words,  for each point $(x,\xi), x\in M,\xi\in T_x^*M$,  we have a linear map  $\sigma(x,\xi):E_x\to F_x$.
\end{definition}

 Assume first that $M$ is a compact manifold.
 To the pseudo-differential operator $A$, one associates  its principal {\it  symbol  $\sigma_p$} which is defined outside the zero section of $T^*M$.
 The operator $A$ is said to be $G$-transversally
elliptic if its principal symbol $\sigma_p(x,\xi)$ is invertible for all $(x,\xi)\in T^*_G M$
such that $\xi\neq 0$.
Using  a $G$- invariant function $\chi$ on $T^*M$ identically
equal to $1$ in a neighborhood of $M$ and compactly supported,
then $\sigma(x,\xi):=(1-\chi(x,\xi)) \sigma_p(x,\xi)$ is  defined on the whole space $T^*M$.
Furthermore  $\sigma(x,\xi)$  restricted to $T_G^*M$  is  an isomorphism
    outside a compact $G$-invariant subset  of   $T^*_GM$.  Thus, by restriction to $T^*_GM$, the symbol
$\sigma$ defines a $K$-theory class $[\sigma]$ in the topological equivariant $K$ -theory group  $K_G^0(T^*_GM)$. This class does not depend of the choice of $\chi$. We still say that this class $\sigma$ is the symbol of $A$.

Let $\hat G$ be the set of equivalence classes of finite dimensional irreducible representations of $G$, and let  $\mathcal C[\hat G]$  be the group  of  $\mathbb Z$-valued functions on $\hat G$.
Let $\chi_\tau(g)={\rm Tr}(\tau(g))$  be the character of the representation $\tau\in \hat G$ of $G$.
We associate to  an element $f\in \mathcal C[\hat G]$  a formal  (virtual) character
$\Theta(f)=\sum_{\tau} f(\tau) \chi_\tau$, that is a formal combination of
the characters $\chi_\tau$  with multiplicities $f(\tau)\in \mathbb Z$.
When $f(\tau)$ satisfies certain moderate growth conditions, then the series
$\Theta(f)(g)=\sum_{\tau} f(\tau) \chi_\tau(g)$ converges, in the distributional sense, to a generalized function on $G$.

The {\it index map} associates to a  transversally elliptic operator  $A$ an element of  $\mathcal C[\hat G]$ constructed as follows.  For every $\tau\in \hat G$, the space $\hom_G(\tau,\ker(A))$ is finite dimensional of dimension $m(\tau,A)$.  Thus $m(\tau,A)$ is the multiplicity of $\tau$ in  the space $\ker(A)$ of smooth solutions of $A$.
We choose a $G$-invariant metric on $M$ and $G$-invariant Hermitian structures on $E,F$. Then
$A^*: \Gamma(M,F)\to \Gamma(M,E)$ is also transversally elliptic.

\begin{definition}
The index  multiplicity  of  the pseudo-differential operator  $A$ is the function $ind_m(A)\in    \mathcal C[\hat G]$ defined by
$$ind_m(A)(\tau):=m(\tau,A) -m(\tau,A^*).$$
\end{definition}

It follows also from Atiyah-Singer \cite{At} that the series  $\sum_{\tau}m(\tau,A) \chi_\tau(g)$ defines a generalized function on $G$.  Thus  we may also  associate to $A$ the generalized function
$$ind(A)(g)=\sum_{\tau}ind_m(A)(\tau)\chi_\tau(g)$$  on $G$ with integral Fourier coefficients.
One of the main points in the  index theory consists in showing that the index factors through the symbols and defines a homomorphism of $R(G)$ modules  from $K_G^0(T_G^*M)$  to  $\mathcal C[\hat G]$.

If $j:U\to M$ is an open $G$-invariant set of a compact $G$ manifold $M$,
we still denote by $j$ the corresponding open embedding from $T^*_GU$ to $T^*_GM$.
 Then $j_{*}$ defines a map from
$K_G^0(T_G^*U)$ to $K_G^0(T_G^*M)$.
 The index of $\sigma\in K_G^0(T^*_G U)$ is defined to be the index of $j_{*}(\sigma)$. The excision property of the index shows that this is independent of the choice of the open embedding $j$ and thus allows  us to define the index map also for  manifolds which can be embedded as open sets of compact ones.
 
 In particular, if $V$ is a vector space with a linear  action of a compact group $G$,
 then $V$ is diffeomorphic to  the  sphere, minus a point.
 Thus we can define the index of any $\sigma\in K_G^0(T_G^*V)$.
More generally, if $U$ is an open $G$-invariant subset of a vector space, we can define the index of $\sigma\in K_G^0(T_G^*U)$.

The problem of computing the index can be reduced, at least theoretically, to the case in which  $G$  is a torus. For a given compact manifold  $M$,  one embeds $M$ into a linear representation and then  is reduced to perform the computations  in the representation.
 \smallskip

In this article, the group $G$ is an abelian compact Lie group.
An irreducible representation $a$ of $G$ is  a one dimensional complex vector space $L_a$, where $G$ acts  via a character  $\chi_a:G\to S^1$, so that $\hat G$ is identified with the abelian group of characters, denoted by $\Gamma$.

\begin{definition}\label{MX}
Let $X$ be a finite list of elements of  $\Gamma$.
Define  the complex vector space   \begin{equation}\label{SMX}
M_X:=\oplus_{a\in  X}L_a.
\end{equation}
\end{definition}

The space $M_X$ is a $G$-manifold and our goal is  the determination of $K_G^0(T^*_GM_X)$.
The basic tool that we shall use is the space of functions $DM(X)$ on $\Gamma$.

\begin{remark}

Let $G$ be a torus  acting on a real vector space $M$ without fixed non zero subspace.
Then $M$ can be given a complex structure, so that the $G$-manifold $M$  is isomorphic  to the space $M_X$ for some list $X$ of weights.
It will be clear that our description of  $K_G^0(T^*_GM_X)$  depends  only of  the list $X$
up to signs.
\end{remark}

\subsection{Examples of transversally elliptic symbols.\label{example}}

Let $M$ be a real vector space provided with a linear action of a compact abelian Lie group $G$, which we may assume to be orthogonal with respect to some chosen Euclidean structure $(v,w)$.
Let $U$ be the Lie algebra of $G$. We assume that there exists $u\in U$ such that the infinitesimal action $\rho(u)$ of $u$ on $M$ is invertible. Such a $u$ will be called regular.  Since  $-\rho(u)^2$ is a symmetric and positive operator,  in particular it is semisimple with  positive eigenvalues  and we can take its unique   square root with  positive eigenvalues.
We choose on   $M$  the  complex structure $J_u=\rho(u)/{(-\rho(u)^2)^{1/2}}$. This complex structure depends only of the connected component $F$ of $u$ in the space of  regular elements.   The given Euclidean structure is the real part of a positive Hermitian structure  $ (v,w)+i(v,J_uw)$  for which the action is unitary.
Let us write $M=\oplus_{a\in X} L_a$, where $X$ is the list of weights of $G$ in the complex vector space $M$.
The connected components of the space of regular elements  are the open faces of the arrangement
$\mathcal H_X$.
By definition of the complex structure, all weights  $a$ are positive  on $u\in U$. We can then define the generalized functions $\Theta_X^{\pm F}(g)$ on $G$ as explained in subsection \ref{convoli}.

$\bullet$ {\bf The tangential Cauchy-Riemann operator.}\label{CR}
With the given Hermitian structure   on $M$,    let $S$ be the unit sphere of $M$. Let  $\mathbb P(M)$ be the complex projective space of $M$.
Consider  on $S$  the  differential operator  $\delta$ acting  on  the  pull back of the Dolbeault complex  on the associated projective space $\mathbb P(M)$ using $\overline \partial + \overline \partial ^*:\sum \Omega^{0,2p}\to \sum  \Omega^{0,2p+1}$.
Then $\delta$ is a $G$-transversally elliptic  differential operator (the tangential Cauchy-Riemann operator) on $S$.
Indeed,   using the Hermitian structure, identify   $T^*S$ with its tangent bundle $TS\subset TM$,   the  subspace $H_p$ of $T^*_pS$ orthogonal to the line $\mathbb  R J_up$ is then identified to the complex subspace of $M$, orthogonal under the Hermitian form to $p$. We call it the horizontal cotangent space.
The symbol of $\delta$ is $\sigma(p,\xi)=c(\xi^1)$ where $\xi^1$ is the projection of $\xi$ on $H_p$, and $c$ the Clifford action of $H_p$ on  $\bigwedge H_p$. This morphism  is invertible if $\xi^1\neq 0$.
We have also  $H_p\oplus  \mathbb R \rho(u)p=T^*_pS$, as the eigenvalues of $-iu$ on $M$ are all positive. Thus we see that $\sigma(p,\xi)$  restricted to $T_G^*S$ is invertible outside the zero section.

The following formula is proven in \cite{At} (Proposition 5.4).
\begin{theorem}
Let $M$ be provided with the complex structure $J_u$ and let $\delta$ be the tangential Cauchy-Riemann operator on the unit sphere of $M$. Then
$$index(\delta)(g)=(-1)^{|X|} g^{a_X}( \Theta_{X}^{-F}(g)-\Theta_X^{F}(g))$$
where $a_X=\sum_{a\in X}a$.
\end{theorem}

\begin{proof}
We recall briefly the proof.
Let $S^1$ be the circle group acting by homotheties on $M$.
We decompose solution spaces with respect to characters $t\mapsto t^n$ of $S^1$.
The  group $G$ acts on
$ \mathbb P(M)$ and on every line bundle $\mathcal O(n)$ on $\mathbb P(M)$.
Thus  the index as an index of $G\times S^1$ is the sum of the index of $G$ in the
cohomology on $\mathbb P(M)$ of the line bundles $\mathcal O(n)$.   Define
$\chi_n(g)$ as the virtual character (as a representation of $G$) in the virtual finite dimensional vector space
$\sum (-1)^i H^{0,i}(\mathbb P(M),\mathcal O(n))$.
Then $$index(\delta)(g)= \sum_{n\in \mathbb Z}\chi_n(g).$$
Let us  show that
\begin{equation}\label{pos}
\sum_{n\geq 0}\chi_n(g)=(-1)^{|X|} g^{a_X}\Theta_{X}^{-F}(g),
\end{equation}
\begin{equation}\label{neg}
\sum_{n<0}\chi_n(g)=(-1)^{|X|+1} g^{a_X}\Theta_{X}^{F}(g).
\end{equation}

For $n\geq 0$, $\mathcal O(n)$   has only 0-cohomology and $H^{0,0}(\mathbb P(M),\mathcal O(n))$ is just the space of  homogeneous polynomials
on  $M$ of degree $n$. So  $\sum_{n=0}^\infty  \chi_n$ is the character of the symmetric algebra  $S[M^*]=\prod_{a\in X}S[L_{-a}]$.
The function $\sum_{k=0}^{\infty} g^{-ka}$ is the character of the action of $G$ in $S[L_{-a}]$.
The function $\Theta_X^{-F}$  is  the product of the functions $-g^{-1}\sum_{k=0}^{\infty}g^{-k a}$.
Thus we obtain Formula (\ref{pos}).
On the other hand, if $n<0$,
we have two cases. If $-|X|-1<n\leq -1$,
then $H^{0,i}(\mathbb P(M),\mathcal O( n))=0 $ for every $i$ . Otherwise we   apply  Serre's duality
  and we obtain the second equality (\ref{neg}).

\end{proof}

$\bullet$\  {\bf Atiyah-Singer pushed symbol.}
We identify $T^*M$ with $M\times M$, using the Hermitian metric on $M$.
Let $c(v):\bigwedge^{even}M\to \bigwedge^{odd}M$ be the Clifford action (\ref{cli}) of $M$ on $\bigwedge M$.
Given, as before,  a regular element $u$  in the Lie algebra of $G$ and letting $\rho(u)$ denote its infinitesimal action on $M$,  we define
\begin{definition}\label{At}
$$At_u(v,\xi)=c(\xi+\rho(u)v).$$
\end{definition}
The morphism $At_u(v,\xi)$ is invertible except if $\xi+\rho(u)v=0$. If furthermore $\xi$ is in $T_G^*M$, $\xi$ is orthogonal to the tangent vector $\rho(u)v$. Thus the support of $At_u(v,\xi)$  restricted to
$T_G^*M$ is the unique point  $v=0,\xi=0$ and   $At_u$ determines an element of $K_G^0(T^*_G M)$, which  depends only of the connected component $F$  of $u$ in the set of regular elements.  We denote it  by $At_F$.
The index of $At_F$ is computed in \cite{At} (Theorem 8.1). In more detail, in the Appendix of \cite{BV1}, it  is constructed an explicit   $G$-transversally elliptic pseudo-differential operator $A$ on the  product of the projective lines  $\mathbb P(L_a\oplus \mathbb C)$. If $j:M_X\to\prod_{a\in X} \mathbb
P(L_a\oplus \mathbb C)$ is the natural open embedding,
  it is shown   that  $j_*(At_F)$ is homotopic  to  the symbol of $A$.
   By definition,  the index of  $At_F$  is that of $A$ and   one has the explicit formula:
\begin{theorem}\label{indexAt}
Let $M$ be provided with the complex structure $J_u$  and let  $At_F\in K_G^0(T^*_GM)$ be the ``pushed'' $\overline \partial$ symbol. Then
$$index(At_F)(g)=(-1)^{|X|}g^{a_X}\Theta_X^F(g).$$
\end{theorem}

\subsection{Some properties of the index map\label{sompr}}
We now recall briefly some  properties of the index map: $K_G^0(T^*_GM)\to C^{-\infty}(G)$ that we will use.
Here again $G$ is an abelian  compact Lie group and $M$ a $G$-manifold, which can be embedded as a.
$G$-invariant open subset of a compact $G$-manifold $\tilde M$.
This hypothesis on $M$  is in place in the rest of the article.
It is stable under the following operations:

 i)\quad If  $U$ is an open $G$-invariant subset of $M$, then $U$ satisfies our hypothesis.

\smallskip

ii)\quad
If $W$  is a real vector space with a linear representation of $G$, then $M\times W$
 satisfies our hypothesis.

\smallskip

iii)\quad
Let $H$ be a closed subgroup of $G$.
Let $M$ be a space with  $H$ action,  open $H$-invariant subset of a compact $H$-manifold.
Let $N=G\times_H M$ be the $G$ space with typical fiber $M$ over
$G/H$.
Then $N$ satisfies our hypothesis.

We now list some properties of the index.

\smallskip

 i)\quad  Any element $\sigma\in K_G^0(T^*_GM)$ arises from  the restriction to $T^*_GM$
of a  $G$- bundle morphism $\sigma(x,\xi):E_x\to F_x$,  such that $supp(\sigma)\cap T^*_G M$ is a compact set. Here  $E,F$ are $G$-equivariant complex vector bundles over $M$.\smallskip

\smallskip

ii)\quad  If $j:\mathcal U\to M$ is an open embedding,  the  map  $j_*: K_G^0(T^*_G\mathcal U)\to  K_G^0(T^*_GM)$ is compatible with the index.

\smallskip

iii)\quad
Let $W$ be  a real vector space with a linear representation of $G$ and $W'$ be the dual vector space. We identify $TW=W\times W$ with $W_{\mathbb C}$ by $(v,w)\to v+iw$. Furthermore, we  identify $TW=W\times W$ with $T^*W=W\times W'$ using an Euclidean structure on $W$. Thus the Bott symbol
$c_{W_{\mathbb C}}(v+i\xi)$  acting on $\bigwedge W_{\mathbb C}$  defines a $G$-equivariant elliptic symbol on $W$.
Its $G$-equivariant index is identically equal to $1$.

Let $i:N\to N\times W$ be the injection of a $G$-manifold  $N$ into
$N\times W$. Then we obtain a map $i_{!}:K_G^0(T^*_GN)\to K_G^0(T^*_G(N\times W))$
given at the level of symbols by $\sigma\mapsto \sigma\odot c_{W_{\mathbb C}}$.
The index of $\sigma$ is equal to the index of $i_{!}\sigma$.

\smallskip

iv)\quad
In case $W=\mathbb R$ with the trivial action,  $T^*_G(N\times\mathbb R)=T^*_GN\times T^*\mathbb R$ and thus $i_!$ is   an isomorphism by Bott periodicity.

\smallskip

v)\quad
Let $H$ be a closed subgroup of $G$.
Then   there is a surjective map
$\hat G\to \hat H$ induced by the restriction of characters.
The dual map induces an injection ${\rm Ind}_H^G: \mathcal C[\hat H]\to \mathcal C[\hat G]$.

Let $M$ be a space with  $H$ action (open subset a compact $H$-manifold), and let $N=G\times_H M$ be the $G$ space with typical fiber $M$ over
$G/H$.  It is easy to see that there is an isomorphism
\begin{equation}\label{indsigma}
i_H^G:K_H^i( T_H^*M)\to K_G^i( T_G^*N)
\end{equation}
and, by
 (\cite{At}, Theorem 4.1), for any $\sigma\in K_H^0(T^*_HM)$,
\begin{equation}\label{indindex}
ind_m(i_H^G(\sigma))={\rm Ind}_H^G (ind_m(\sigma)).
\end{equation}

We shall also need the following simple consequence of the previous facts:

\begin{lemma}\label{genrem}Let $G$ be a compact Lie group and $\chi:G\to S^1$ be a surjective character.
Set $H:=\ker  \, \chi$ be  the kernel of $\chi$.

 Take a manifold $M$  over which $G$ acts and consider the product $\mathbb C^*\times M$, with  the action of $G$ on the first factor induced by  $\chi$.

There is    an  isomorphism
\begin{equation}\label{prii1}
k: K^{i }_H(T^*_HM)\cong K^i_G(T^*_G(\mathbb C^*\times M)) .
\end{equation}
Moreover,  if $\sigma\in K^{0 }_H(T^*_HM)$, we have $ind_m(k(\sigma))={\rm Ind}_H^G (ind_m(\sigma)).$\end{lemma}
\begin{proof} Since $\mathbb C^*=S^1\times \mathbb R^+$, we get
by iv) that the inclusion $i:S^1\times M\to \mathbb C^*\times M$ induces  the isomorphism $$i_!:K_G^i(T^*_G(S^1\times M))\to K_G^i(T^*_G(\mathbb C^*\times M))$$ which at the level of $K^0$ is  compatible with the index.

On the other hand,  the space  $G\times_H M$ identifies with  $S^1\times M$  via the map
  $[g,m]\mapsto [\chi(g),g\cdot m]$. So  (\ref{indsigma}) gives us the isomorphism
   $$i_H^G: K_H^i(T^*_H(M)) \to K_G^i(T^*_G(S^1\times M)).$$ If $\sigma\in K_H^0 (T^*_HM)$,   then $ind_m(i_H^G(\sigma))={\rm Ind}_H^G (ind_m(\sigma))$ by
   Formula (\ref{indindex})  .

Thus we can take  $k:=i_!i_H^G$.
\end{proof}

\section{ Equivariant $K$-theory and  Dahmen-Micchelli spaces}
This section contains the main results of this paper, that is Theorem \ref{first}     and Theorem \ref{ilprinc}.
\subsection{Two exact sequences\label{tes}}
Let $G$ be, as before, a compact abelian Lie group
of dimension $s$ and
$M_X:=\oplus_{a\in  X}L_a$ as in (\ref{SMX}).
We  assume that $X$ is a non degenerate list of characters of $G$.

Given a vector $v\in M_X$,  its {\it support} is the sublist  of elements $a\in X$ such that $v$ has a non zero coordinate in   the
summand $L_a$.

 If $Y$ is   the support of  $v$,  an element $t$ of  $G$  stabilizes $v$ if and only if  $t^a=1$ for all $a\in Y$.  If $Y$ spans a  rational subspace of dimension $k$,   the $G-$orbit of $v$ has dimension $ k$.

For any rational subspace $ \underline r$,
we may consider the subspace  $M_{\underline r}:=\oplus_{a\in X\cap \underline r}L_a$ of $M_X$. We set
\begin{equation}
M_{\leq i}:=\cup_{\underline r\in S_X^{(i)}}M_{\underline r},\quad M_{\geq i}:=M_X\setminus M_{\leq i-1},$$ $$ F_i :=M_{\leq i}\setminus M_{\leq i-1}=M_{\geq i}\setminus M_{\geq i+1}=M_{\leq i}\cap M_{\geq i}.
\end{equation}   Notice that   $$  M_X= M_{\geq 0}\supset M_{\geq 1}\supset M_{\geq 2}\supset\cdots \supset M_{\geq s} :=M^f_X.$$
 The set $M_{\leq i}$ is the closed set of points in $M$ with the property that the orbit has dimension $\leq i$ while $M_{\geq i}$ is the open set of points in $M$ with the property that the orbit has dimension $\geq i$. The set   $F_i$ is open in   $M_{\leq i}$ and closed   in $M_{\geq i}$ and it is  the set of points in $M$ whose orbit under $G$ has dimension exactly $ i$.
In particular,  $F_s=M_{\geq s}=M^f_X$ is the open set of points in $M$ with finite stabilizer  under the action of $G$, which plays  a particular role.

Given a rational subspace $\underline r$,  we have denoted by   $G_{\underline r}$   the subgroup of $G$ joint kernel of the   elements in  $\Gamma\cap\underline r$. The group $ G_{\underline r}$ is a torus and acts trivially on $M_{\underline r} :=\oplus_{a\in X\cap \underline r}L_a$ inducing an action  of $G/G_{\underline r}$.
\begin{definition} We define the set $F(\underline r)$  to be the open set of $M_{\underline r}$ where $G/G_{\underline r}$ acts with  finite stabilizers.

In other words, the connected component  of the stabilizer of an element of $F(\underline r)$
is exactly the group $G_{\underline r}$.\end{definition}\begin{remark}  By definition of $G_{\underline r}$, the set  $F(\underline r)$ is non--empty.
The set $F_i$ is the disjoint union of the  sets  $F(\underline r)$ as $\underline r$ runs over all rational subspaces of  dimension  $i$.
Thus the space  $M_X$ is the disjoint union of the locally closed strata  $F(\underline r)$.
\end{remark}

We now analyze the equivariant $K$  theory of $T^*_G M_X^f$.

Let  $a\in X$ be an element of   infinite order so that  the homomorphism  $g^a:G\to S^1$ is  surjective.
Set $Z:= X\setminus\{a\}$ and $G_a:=\ker   \, g^a$.
 Denote by $\tilde   Z$  the list of the restrictions to  $G_a$ of the elements of $Z$.
    For $v\in M_X$,
  denote by $v_a\in \mathbb C$ its coordinate in $L_a$ with respect to a choice of a basis of the one dimensional vector space $L_a$.

   The set   $M_Z^f:=\{v\in M^f_X\,|\, v_a=0\}$ is closed in $M_X^f$.  Denote by $i: M_Z^f\to M_X^f$ the closed embedding and by $j:M_X^f\setminus M_Z^f \to M^f_X$ the open embedding of the complement.

\begin{lemma}\label{kiso}
There exists an isomorphism
$$k:K_{G_a}^i(T_{G_a}^*M_{\tilde Z}^f)\to   K_G^i(T_G^*(M_X^f\setminus M_Z^f)).$$
If $\sigma\in K_{G_a}^0(T_{G_a}^*M_{\tilde Z}^f)$, we have $ind_m(k(\sigma))={\rm Ind}_{G_a}^G (ind_m(\sigma)).$\end{lemma}
\begin{proof}
Take an element  $(v_a,w)\in L_a\times M_Z$  with $v_a\neq 0$, its stabilizer  in $G$ is the subgroup of $G_a$ stabilizing $w$, therefore
the space
$M_X^f\setminus M_Z^f$ is  isomorphic to $\mathbb C^*\times M_{\tilde Z}^f$.
Thus  we are  in the setting of Lemma   \ref{genrem}.
\end{proof}
For a real vector space $W$,  we shall denote by $W'$ its dual.
Consider the projection
$p:T_G^*M_X^f\to M_X^f$.
Then $p^{-1} M_Z^f$ is a closed subset of $T_G^* M_X^f\subset M_X^f\times M_X'$
and $T_G^*M_X^f\setminus p^{-1}M_X^f$ is equal to $T_G^*(M_X^f\setminus M_Z^f)$.
We use the same notations $i,j$ also in this setting for the closed and open embedding associated.
Remark the following fact.
\begin{lemma}\label{kiso1}
\begin{enumerate}
\item
We have $p^{-1} M_Z^f= T_G^*M_Z^f\times L_{a}'$.

\item

We have an isomorphism $C_a:K_G^i(T^*_GM_Z^f)\to K_G^i(p^{-1}M_Z^f)$. \end{enumerate}
\end{lemma}
\begin{proof}
The first assertion is immediate to verify. The second follows from the first  and Thom  isomorphism.
\end{proof}

The first theorem is:

\begin{theorem}\label{ipfT}
\begin{enumerate}\item $K_G^{1}(T^*_GM_X^f)=0$.
\item
If  $a\in X$ has infinite order, there  is a short exact sequence:
\begin{equation}\label{prex1} 0\to K^0_{G_a}(T_{G_a}^*M_{\tilde Z}^f) \stackrel{j_*k}\to  K^0_G(T_G^* M_X^f )\stackrel {C_a^{-1}i^*} \to K^0_G(T_G^* M_Z^f )\to 0. \end{equation}\end{enumerate}
\end{theorem}
\begin{proof}
If $G$ is finite,  then $M^f_X=M_X$ and the first statement is  Bott periodicity  while the second statement does not exist.
So we assume that $G$ has positive dimension.
Therefore we can choose $a\in X$ of infinite order.  By induction we can assume that   $K_G^{1}(T^*_GM_Z^f)=0  =K_{G_a}^{1}(T^*_GM_{\tilde Z}^f)$.

  If $Z$ is degenerate, then $M_Z^f$ is empty,  $j$ is the identity and our claims reduce to Lemma \ref{kiso}.  Otherwise  both statements follow from Lemma \ref{kiso}, Lemma \ref{kiso1} and the long exact sequence of $K-$theory.\end{proof}

We choose an Hermitian product $\langle v,w\rangle$ on $M_X$ and we identify the (real) vector space  $M_X'$,  dual  to $M_X$, with $M_X$ using the real part of the Hermitian product.
Thus  $T^*M_X^f$ is identified with $M_X^f\times M_X'$  , where $M_X'$ is the space $M_X$ with the opposite complex structure, so that  $T^*M_X^f$ is a complex vector bundle over $M_X^f$.
The following proposition allows us to reduce the computation of the equivariant $K$-theory of $M_X^f$ to that of $T_G^*M_X^f$.

\begin{proposition}\label{vanish}

\begin{enumerate}\item $K_G^{s+1}(M_X^f)=0$.
\item
There is a natural isomorphism:
$K_G^s(M_X^f)\to K_G^{0}(T^*_G M_X^f)$.
\end{enumerate}
\end{proposition}
\begin{proof}
Let $M$ be a  manifold with an action   of  $G$.
Assume $G$ has finite stabilizers on $M$. We claim that, for every $i$, there is a natural isomorphism between
$K^{i+s}_G(T_G^* M )$ and  $K^{i}_G( T^* M).$
 In fact the infinitesimal action of ${\rm Lie}(G)$  determines a trivial vector bundle $L$ in $TM$.
  Using a $G$-invariant Riemannian structure on $M$,
  we   identify $T^*M$ with $TM$ so that $T^*_GM$ is identified to  the orthogonal of $L$. Thus  we have the product decomposition
$T_G^*M \times {\rm Lie}(G)=TM$.  In this decomposition $G$ acts trivially on  the $s-$dimensional factor  ${\rm Lie}(G)$. We can apply Bott periodicity
$$K_G^0T_G^*M  =K_G^sTM,\quad K_G^1T_G^*M  =K_G^{s+1}TM.$$
In our case,  $T^*M_X^f$ is a complex vector bundle over $M_X^f$, so  we can apply the Thom isomorphism for this bundle  and   we obtain, using Theorem \ref{ipfT},
that $K_G^s(M_X^f)$ is isomorphic to $K_G^{0}(T^*_G M_X^f)$
and that $K_G^{s+1}(  M_X^f)$ vanishes.
\end{proof}

For a rational subspace $\underline r$,  the action of $G$ on $F(\underline r)$  factors through  $G/G_{\underline r} $ and, with respect to this action,   $F(\underline r)=M_{X\cap \underline r}^f $. Thus $$K^i_G (F(\underline r))=R(G)\otimes_{R(G/G_{\underline r})}K^{i}_{ G/G_{\underline r}} (F(\underline r)).$$
In particular, by    Theorem \ref{ipfT}, we deduce that $K_G^1( T^*_G F(\underline r))=0$.

Now set $\tilde T^*_G F(\underline r):=T_G^*M_X| F(\underline r)$,  the restriction of  $T_G^*M_X$ to $F(\underline r)$.
We see that  $\tilde T^*_G F(\underline r)= T^*_G F(\underline r)\times M_{X\setminus\underline r}'$,
  so we have a Thom  isomorphism
   $$C_{\underline r}:  K_G^0(  T^*_{ G }  F(\underline r))\to K_G^0(\tilde T^*_G F(\underline r)),\quad K_G^1(\tilde T^*_G F(\underline r))=0.$$
Choose $0\leq i\leq s$.
We pass now to study the $G$-invariant open subspace  $M_{\geq i}$ of $M$.
The set  $M_{\geq i+1}$ is open in $M_{\geq i}$ with complement the  space $F_i$ disjoint union of the spaces $F(\underline r)$ with  $\underline r \in S_X^{(i)}$. Denote by $\tilde T_G^* F_i$ the restriction of $T^*_G M$ to $F_i$,  disjoint union of the spaces
$\tilde T^*_G F(\underline r).$
Denote $j:  M_{\geq i+1}\to M_{\geq i}$ the open inclusion and $e: \tilde T_G^* F_i\to T_G^*M_{\geq i}$ the closed embedding. Let $C_i$ be the Thom isomorphism from $K_G^0( T^*_G F_i)$ to $K_G^0(\tilde T^*_G F_i)$ direct sum of the Thom isomorphisms $C_{\underline r}$.

\begin{theorem}\label{Tmenos} For each $0\leq i\leq s-1$,
\begin{enumerate}\item $K_G^{1}(T^*_GM_{\geq i})=0$.
\item
The following sequence is exact
\begin{equation}\label{lasecondina}0\to K_G^{ 0}(T_G^*M_{\geq i+1})\stackrel{j_*}\to K_G^{ 0}(T_G^*M_{\geq i})\stackrel{C_{i}^{-1}e^*}\to K_G^{ 0}(T_G^*F_i)\to 0 .\end{equation}
\end{enumerate}
\end{theorem}
\begin{proof}
Since $M_{\geq s}= M^f_X$, we can assume by induction on $s-i$ that i) holds for each $j>i$. Also by Theorem \ref {ipfT} i)
we get that $K_G^{ 1}(T_G^*F_i)=0$ for each $0\leq i\leq s-1$. Using this both statements follow immediately from the long exact sequence of equivariant K-theory.
\end{proof}

\begin{remark}
The fact that the sequence (\ref {lasecondina}) is exact is
proved in \cite{At} using a splitting. We will comment on this point  in  Section \ref{genK}.
\end{remark}

\subsection{Two commutative diagrams}

Let $N$ be  a complex  representation space for $G$.
Recalling  the structure of $R(G)$ module of the  equivariant $K$-theory, the multiplication
 by the difference $\bigwedge^{even} N-\bigwedge^{odd} N$  will be denoted by $\bigwedge_{-1}N\otimes  -$.
 This is by definition the action of  the element $\det_N(1-g)\in R(G)$ on the equivariant $K$-theory.

\begin{lemma}\label{ico}
Take a sublist $Y$ in $X$ and decompose $M_X=M_{X\setminus Y}\oplus M_Y$.
Let $U$ be an open $G$-invariant set contained in $M_{X\setminus Y}\times ( M_Y\setminus \{0\})$. Then if $\sigma\in K^i_G(U)$  or $\sigma\in K^i_G(T_G^*U)$, we have $\bigwedge_{-1}M_Y \otimes  \sigma=0$.
\end{lemma}
\begin{proof}
We give the proof for  $U$, the case of  $T_G^*U$ being identical.

Take $v\in U$ and decompose it as
$v=v_{X\setminus Y}+ v_Y$ with $v_{X\setminus Y}\in M_{X\setminus Y}$ and $v_Y\in M_Y$.
The component $v_{Y}$ is not zero  by assumption.
Consider the complex  $G$-equivariant vector bundle $V_Y=U\times M_Y$ on $U$. Set now  $ E^+:=U\times \bigwedge^{even} M_Y$, $  E^- =U\times \bigwedge^{odd}M_Y$.
Choosing an Hermitian metric on $M_Y$,  for every $u\in M_Y$, we get the  Clifford action
$c(u): \bigwedge^{even} M_Y\to  \bigwedge^{odd} M_Y$ of $M_Y$ on $\bigwedge M_Y$,
which is an isomorphism  as soon as $u\neq 0$.

Going back to our bundles $E^+, E^-$,  for every $\epsilon\in [0,1]$, define the bundle map
${\bf c_\epsilon}:E^+\to E^-$   by
$${\bf c}_\epsilon(v,\omega)=(v, \epsilon c(v_Y)\omega).$$

If $\sigma \in K_G^0(U)$, the element $\bigwedge_{-1}M_Y\otimes   \sigma\in K_G^0(U)$
is represented by the morphism $ {\bf c_0}\odot  \sigma$, homotopic to
${\bf c_1}\odot \sigma$. This last bundle map is an isomorphism since $v_Y\neq 0$ on $U$.
This implies that
$\bigwedge_{-1}M_Y \otimes  \sigma=0$.
\end{proof}

We apply this to the open set $M_X^f$ where $G$ acts with finite stabilizers.
 If $Y$ is a cocircuit,   $M_{X\setminus Y}\cap M_X^f=\emptyset$,
   thus  $\bigwedge_{-1}M_Y \otimes  \sigma=0$ for all $\sigma\in K^0(T_G^*M_X^f)$.
   As the index map is a $R(G)$ module map, this implies that for  cocircuit $Y$,
   the generalized function $index(\sigma)(g)$ on $G$ satisfies the equation
   $\prod_{a\in Y}(1-g^a) index(\sigma)(g)=0$. The function $ind_m(\sigma)$ on $\hat G$ is
   the Fourier transform of the function $index(\sigma)$. It follows that
   $\nabla_Y ind_m(\sigma)=0$.

    Thus we obtain

\begin{corollary}\label{valind}
The  multiplicity index map $ind_m$ maps  $K^0_G(T_G^*M_X^f) $ to the space $DM(X)$.
\end{corollary}

More generally the same argument shows that

\begin{corollary}\label{laquetio}
Choose  $0\leq i\leq s$.
 If $\sigma\in K_G^0(T_G^*M_{\geq i})$ and $\underline t$ is a rational subspace of dimension strictly less than $i$,
 then  $\bigwedge_{-1}M_{X\setminus \underline t}\otimes  \sigma=0$.
\end{corollary}

Let us now split $X=A\cup B$ and $M_X=M_A\oplus M_B$. Let $p:T^*_G M_X\to M_X$ be the projection and consider   $ \tilde T_G^*M_A:=p^{-1}M_A$.
We  have  $ \tilde T_G^*M_A= T_G^*M_A\times  M_B'$. In particular,  we get a Thom  isomorphism $$C_{M'_B}: K_G^0(T_G ^* M_A)\to K_G^0(\tilde T_G^*M_A)\cong  K_G^0(T_G^*M_A\times  M_B') .$$

Denote by $i$ the closed inclusion $M_A\to M_X$, and, by abuse of notation,  also the inclusion
$\tilde T_G^*M_A\to  T_G^*M_X $ above $i$.
Then $i$ induces   the  morphisms $i^*:    K_G^0(T^*_GM_X )\to K_G^0(\tilde T^*_GM_A)$
 and  $i_!:K_G^0(T^*_G  M_A)\to K_G^0(T^*_GM_X ).$
Combining these 3 maps,   we claim that
\begin{lemma}\label{IMC} Take $\sigma\in K_G^0(T_G^*M_X )$,
then  $i_!C_{M'_B}^{-1} i^*(\sigma)=\bigwedge_{-1} M_B\otimes \sigma$.

\end{lemma}
\begin{proof} Since we are working on vector spaces,
we can assume that all vector bundles are topologically trivial. Thus we can represent $\sigma$ as given by a variable linear map
$\sigma(v,w,\xi,\eta):E\to F$ where $E,F$ are complex  representation  spaces,
  $v\in M_A, w\in M_B,\xi\in M_A',\eta\in  M_B'$.
   Now  $\sigma$ restricts to an element $\tilde \sigma:=i^* \sigma$ in $K^0_G(T^*_GM_X|M_A )$
    which is   represented by the map  $\sigma(v,0,\xi,\eta).$  Since  $T^*_GM_X|M_A =T^*_G M_A \times  M_B' $, the element
     $\tilde \sigma$   is equivalent  to  ${\bf c}_ {M_B'}\odot q^*\tau$, where  $q:T^*_GM_A \times M_B' \to T^*_G M_A $ is the projection, $\tau$  a  transversally elliptic symbol on  $M_A $
     and ${\bf c}_ {M_B'}$ the Bott symbol with support the zero section of the bundle $T^*_GM_A \times M_B'$ on $T^*_GM_A $.
        Thus we have to show that $\bigwedge_{-1}M_B\otimes \sigma$
        and  $i_! (\tau)$   are homotopic.

By definition,  a representative of the symbol $i_!(\tau)$  on $M_A\times M_B$ is  the product of the symbol
 ${\bf c}_{M_B\otimes_{\mathbb R}\mathbb C }$ by the symbol $q^*\tau$.
  As  ${\bf c}_{M_B\otimes_{\mathbb R}\mathbb C }={\bf c}_{M_B}\odot {\bf c}_{M'_B}$,  we see that
$i_!(\tau)=q^*\tau\odot {\bf c}_{M_B}\odot {\bf c}_{M'_B}= \tilde \sigma \odot {\bf c}_{M_B}$.

As we have seen before, the symbol defined as  $\bigwedge_{-1}M_B\otimes \sigma$  on the manifold $M_A\times M_B$ is homotopic to
the element ${\bf c}_{M_B}\odot \sigma$.
Now consider the symbol $\sigma(t) ( v,w,\xi,\eta)=\sigma(v,tw,\xi,\eta)$ on $T^*M$. The intersection of the support of   ${\bf c}_{M_B}\odot \sigma(t)$ with $T^*_G M$  stays compactly supported for all $t$. Indeed its support remains constant: this is the intersection of the support of $\sigma$ with $T^*_G M_A$.
So  we obtain
the desired homotopy between $\bigwedge_{-1}M_B\otimes \sigma={\bf c}_{M_B}\odot \sigma(1)$
        and  $i_! (\tau)={\bf c}_{M_B}\odot \sigma(0)$ and the claim follows.

\end{proof}

With the previous notations, $X=A\cup B$, $i:T_G^*M_A\times  M_B' \to  T_G^*M_X $.
\begin{corollary}\label{nablaindex}
Take
$\sigma\in K_G^0(T_G^*M_X )$.
Let $\sigma_0=C_{M'_B}^{-1} i^*(\sigma)\in K_G^0(T^*_G M_A)$.
Then, we have the equality of generalized functions on $G$:
 $$det_{M_B}(1-g) index(\sigma)(g)=index(\sigma_0)(g).$$
\end{corollary}

We are now ready to compare the exact sequence (\ref{prex1}) with the exact sequence given by Theorem (\ref{lasf}) using the index. We get
\begin{theorem}\label{first} The diagram
\begin{equation}\label{dia1}\begin{CD} 0\to K^0_{G_a}(T^*_{G_a}M_{\tilde Z}^f) @>{ j_* k}>> K^0_G(T_G^* M_X^f ) @>{C_a^{-1}i^*} >>K^0_G(T_G ^*M_Z^f )\to 0\\  @Vind_mVV@Vind_mVV@Vind_mVV@.\\
\hskip-0.6cm 0\longrightarrow DM(\tilde Z ) @>i_a>> DM(X ) @>\nabla_a >>\hskip0.7cm DM( Z)\to 0 \end{CD} \end{equation}
is commutative. Its vertical arrows are isomorphisms.

In particular, the index  multiplicity map gives an isomorphism between $K^0_G(T^*_G M_X^f )$ and $DM(X)$.
\end{theorem}
\begin{proof}
We start by remarking that, by Corollary \ref{valind},  all the vertical maps in our diagram
are indeed taking values in the corresponding Dahmen-Micchelli spaces.

So we need to show commutativity.
To prove the commutativity  of the square on the right hand side, using Fourier transform, we need to prove that
$(1-g^a) index(\sigma)(g)=index(C_a^{-1} i^*(\sigma))(g)$.

From the symbol $\sigma$ on $M_X^f$, an open set in $M_X$, we deduce a symbol on $M_X$  with same index, by the excision property of the index.
Thus  the commutativity follows from Corollary  \ref{nablaindex} applied to $A=Z,\, B=\{a\}$.
As for the square on the left hand side, since $j_*$  is an open embedding, it preserves indices.  The statement thus follows from  Proposition \ref{kiso}.

By induction we can then assume that the two external vertical arrows  are isomorphism so, by the five Lemma, also the central one is and everything follows.
\end{proof}

Summarizing we have isomorphisms$$\boxed{K^{s+1}_G(  M_X^f ) \cong K^1_G(T_G^* M_X^f )=0,\ \ \ \
K^{s}_G(  M_X^f ) \cong K^0_G(T_G^* M_X^f )\cong  DM(X ).}$$

Let us make two obvious remarks on the isomorphism  $K^0_G(T_G^* M_X^f )\cong  DM(X ).$

\begin{remark}
The space $K_G^0(T^*_G M_X^f)$ depends only of the manifold $M_X$ considered as a real manifold.
 By Remark \ref{real}, the space $DM(X)$ depends only of the list $X$ up to change of signs.
\end{remark}

\begin{remark}\label{supportfourier}
If $U$ is a $G$ manifold, the index of an element $\sigma\in K_G^0(T^*_G U)$ is a generalized function  supported on the set of points $g\in G$ such that $g$ has a fixed point in $U$.

We have seen in \S \ref{supports}  that Fourier transforms of elements in $DM(X)$ are supported on the finite set of points $P(X)$. This is in agreement with the fixed point philosophy that we just recalled.  In fact,
an element  $g\in G$ has a fixed point $v$ in $M_X^f$ if and only if $g\in P(X)$. Indeed if   $g\in P(X)$, there exists a basis $\underline b$ of $V$ extracted form $X$  with $g^{b_i}=1$, for all $b_i\in \underline b$. Thus any element
$v\in M_X$ with non zero coordinates on each $L_{b_i}$ is fixed by $g$, and is in $M_X^f$.
\end{remark}

We now come to our next commutative diagram.
\begin{lemma}\label{valindTotal} For each $s\geq i\geq 0$, the index multiplicity map $ind_m$ sends
$K_G^0(T^*_G M_{\geq i})$ to the space
$\tilde {\mathcal F}_{i}(X)$.
\end{lemma}
\begin{proof}

Recall that  $\tilde {\mathcal F} _{i}(X)$ is the subspace in $\tilde {\mathcal F}(X)$ such that
$\nabla_{X\setminus \underline t}f=0$ for all $\underline t\in S_X^{(i-1)}$. Denote  by $\ell:\tilde {\mathcal F} _{i}(X)\to \tilde {\mathcal F}(X)$ the inclusion.

By Corollary \ref{laquetio}, if $\sigma\in K_G^0(T_G^*M_{\geq i})$ and $\underline t$ is  a rational subspace of dimension strictly less than $i$,  we have $\bigwedge_{-1}M_{X\setminus \underline t}\otimes  \sigma=0$. Thus  $\nabla_{X\setminus \underline t}ind_m(\sigma)=0$.
It follows that the only thing we have to show is that, if $\sigma\in K_G^0(T^*_GM_X)$,  then
$ind_m(\sigma)$  lies in $\tilde{\mathcal F}(X)$.
Take a rational subspace $\underline r$.  By Lemma  \ref{IMC}, the index of $\bigwedge_{-1}M_{X\setminus \underline r}\otimes  \sigma$ equals the index of an element $\sigma_0\in K_G^0(T^*_GM_{X\cap\underline r})$. But the action of $G$ on $M_{X\cap\underline r}$ factors though the quotient $G/G_{\underline r}$
whose character group is $\Gamma_{\underline r}$.
Thus $K_G^0(T_G^*M_{X\cap\underline r})\cong R(G)\otimes_{R(G/G_{\underline r})}K^0_{G/G_{\underline r}}(T^*_{G/G_{\underline r}}M_{X\cap\underline r})$,
hence  $\nabla_{X\setminus \underline r}ind_m(\sigma)=
ind_m(\sigma_0)$ lies in $R(G)\otimes_{R(G/G_{\underline r})}\mathcal C[\Gamma_{\underline r}]$ as desired.
\end{proof}

Our second commutative diagram and  main  theorem characterizes the values of the index on  the entire $M_X$. This time,  we use the notations and the exact sequences contained in Theorem \ref{Tmenos} and Corollary \ref{nablaindex}
\begin{theorem}\label{ilprinc} For each $0\leq i\leq s$,
\begin{itemize}
\item the diagram
$$\begin{CD}  0\to K_G^0(T_G^*M_{\geq i+1})@>{j_*}>> K_G^0(T_G^*M_{\geq i}) @>{C_{i}^{-1}e^*}>>K_G^0(T_G^* F_i)\to 0\\ @Vind_mVV@Vind_mVV@Vind_mVV@.\\
\hskip-1.0cm 0\to\tilde{\mathcal F}_{i+1}(X) @>\ell >> \tilde{\mathcal F}_{i}(X) @>\mu_i>>
\oplus_{\underline r\in S^{(i)}_X}DM^{(G)}(X\cap
\underline r )\to 0 \end{CD} $$
commutes.
\item Its vertical arrows are isomorphisms.

\item In particular, the index gives an isomorphism between $K_G^0(T^*_GM_X)$ and $\tilde {\mathcal F}(X)$.

\end{itemize}

\end{theorem}
\begin{proof} Lemma \ref{valindTotal} tells us that the diagram is well defined. We need to prove commutativity.

Again, we prove that  the square on the right hand side is commutative using   Corollary \ref{nablaindex}.
The square on the left hand side is commutative since $j_*$ is compatible with the index and $\ell$ is the inclusion.

Recall that $K_G^0(T_G^*M_{X\cap\underline r})\cong R(G)\otimes_{R(G/G_{\underline r})}K_{G/G_{\underline r}}^0(T^*_{G/G_{\underline r}}M_{X\cap\underline r})$ and that $DM^{(G)}(X\cap
\underline r )\cong R(G)\otimes_{R(G/G_{\underline r})}DM(X\cap \underline r)$. Using Theorem \ref{first},  this implies that the right vertical arrow is always an isomorphism.

We want to apply descending induction  on $i$. When $i+1=s$, since $M_{\geq s}= M_X^f$ and $\tilde{\mathcal F}_{s-1}(X)=DM(X)$, Theorem \ref{first}   gives that the left vertical arrow is an isomorphism. So assume that  the left vertical arrow is an isomorphism. We then deduce by the five Lemma that the central vertical arrow is an isomorphism and conclude by induction.
\end{proof}

\begin{remark}   Fourier transforms of elements in
 $\tilde{\mathcal F}_{i}(X)$  are supported on the closed set  $P_{s-i}(X)$ described in
 Proposition \ref{supportFi}.
  This is again in agreement with the fixed point philosophy that we  recalled in Remark \ref{supportfourier}.  In fact,
an element  $g\in G$ has a fixed point $v$ in $M_{\geq i}$ if and only if $g\in P_{s-i}(X)$.
\end{remark}

\section{ Generators of $K_G^0(T^*_GM_X)$}\label{genK}

In this section, we show that the generators of $\tilde{\mathcal F}(X)$ constructed in Section \ref{genf} corresponds via the index map to the generators of $K_G^0(T^*_G M_X)$ constructed by Atiyah-Singer in \cite{At}.

We assume that $X$ does  not contain any element of finite order. This is harmless as otherwise we need only to tensor our results with the Bott symbol on $\oplus_{a\in \Gamma_t}L_a$.

Recall that, given a   manifold $M$,   a way to construct   elements of $K_G^0(T^*_GM)$ is to take a closed  $G-$ manifold $N$ embedded  by   $i:N\to M$. Then we have a  map $i_!:K_G^0(T^*_GN)\to K_G^0(T_G^*M)$.

For our case $M=M_X$, we shall take the following   manifolds.
Take a flag  $\phi$ of rational subspaces $0=\underline r_0\subset \underline r_1\subset
 \underline r_2\subset \cdots  \subset \underline r_s$ with $\dim(\underline r_i)=i$ (and $s=\dim G$).
 Consider then the spaces $E_i:=\oplus_{a\in (X\cap \underline  r_i)\setminus \underline r_{i-1}}L_a $.
We choose an orientation for each $\underline r_i$  and divide  the set of characters  $Z_i:=(X\cap \underline  r_i)\setminus \underline r_{i-1} $ into positive and negative elements  $A_i,B_i$.   Accordingly,  we change the complex structure  on each $L_b$ for which $b$ is negative into its conjugate  structure. Let $A$  be the union of the sets $A_i$ and $B$ the union of the sets $B_i$.
Choose a $G$-invariant Hermitian metric $h_i$ on $E_i$ and consider the unit sphere $S_i(h_i)$ on $E_i$.
The product $S_\phi(h)=\prod_{i=1}^s S_i(h_i)$ is a closed submanifold  of $M_X^f$.

The tangent space at a point $p$  of $S_\phi(h)$ decomposes as a vertical space generated by the rotations on each factor $E_i$ and  the horizontal space $H_p$, a lift of the tangent space to the corresponding product of projective spaces. The horizontal space is a Hermitian vector space. The tangential Cauchy-Riemann operator  $\delta_{\phi}$  is a differential operator on $S_\phi(h)$,   the product of the operators $\delta_i$ described in Subsection \ref{CR}. The index
of $\delta_{\phi}$ is the product of the indices  of $\delta_i$.

Let $c_\phi(p,\xi)= c(\xi^1)$
be the  Clifford action on $\bigwedge H_p$  of the projection  $\xi^1$ of $\xi$ on  the horizontal tangent space $H_p$.

We then have \cite{At} the following theorem.
\begin{theorem}
$$\,index(\delta_{\phi})(g)= (-1)^s(-1)^{|B|} g^{\sum_{a\in A}a}
\theta_{\phi}^X(g).$$
\end{theorem}
In fact  $c_\phi=c_1\odot c_2\odot \cdots \odot c_s$ is the external   product of the symbols $c_i$ of the operators
$\delta_i$.

Let $i_\phi$ be the closed embedding of $S_\phi(h)$ in $M_X^f$.
We can then give an  ``easy" proof of the following theorem of Atiyah-Singer (Theorem 7.9  of \cite{At})
\begin{theorem}
Let $i_\phi$ be the closed injection of $S_\phi(h)$ to $M_X$.
Then the  elements $(i_\phi)_{!} c_\phi$ generate  $K_G^0(T^*_G M_X^f)$.
\end{theorem}
\begin{proof}
This follows from Theorem \ref{LABAS}  giving  generators for   $DM(X)$, and the fact that the index multiplicity  map is an isomorphism onto $DM(X)$.
\end{proof}
\begin{remark}
In fact   Theorem \ref{LABAS}  gives a basis of  the space $DM(X)$,  so that the previous theorem can be refined accordingly.
\end{remark}

According to \S \ref{example}, if we consider $M_X$ as a real representation of $G$, each connected component $F$ of the space of regular elements $u\in U=Lie(G)$ gives us a complex structure $J_F$ on $M_X$ and a corresponding
``pushed" symbol $At_F$ with  $index(At_F)(g)=(-1)^{|B|}g^{\sum_{a\in A} a} \Theta_X^F $  (Theorem \ref{indexAt}).
\begin{theorem}\label{crucial}
The symbols $At_F\in K_G^0(T^*_GM_X)$, where $F$ varies over all
open faces of the arrangement $\mathcal H_X$,  give  us a set of generators for
$K_G^0(T^*_G M_X)$
\end{theorem}
\begin{proof}
This  follows from Theorem \ref{GENF}  giving  generators for the space   $\tilde{\mathcal F}(X)$, and the fact that the index multiplicity map is an isomorphism onto $\tilde{\mathcal F}(X)$.
\end{proof}

\begin{remark}
After checking naturality  axioms,  Theorem \ref{crucial} reduces the proof of the cohomological index formula given by
\cite{BV1}
or \cite{parver} to the case of the symbols $At_F$. So  Theorem \ref{crucial}  is  crucial
in establishing a cohomological formula valid for any transversally elliptic operator.
\end{remark}

\begin{remark}
Consider the exact sequence  in Theorem \ref{Tmenos}.
$$0\to K_G^{ 0}(T_G^*M_{\geq i+1})\stackrel{j_*}\to K_G^{ 0}(T_G^*M_{\geq i})\stackrel{C_{i}^{-1} e^*}\to K_G^{ 0}(T_G^*F_i)\to 0 .$$

In \cite{At}, the exactness of this sequence
is  proved using a splitting. We recall the proof  of \cite{At} for $i=s-1$, the proof for any $0\leq i\leq s$ being   identical.
Consider $\underline r$ of codimension $1$  and write $M_X=M_{X\cap \underline r}
\oplus M_{X\setminus \underline r}$.
Choose a regular element $u_{\underline r}$ vanishing on $\underline r$.
We modify the complex structure on $M_{X\setminus \underline r}$ so that all the weights of $G$  on $M_{X\setminus \underline r}$ are positive on $u_{\underline r}$. We then construct the corresponding ``pushed" $\overline \partial$ symbol on $M_{X\setminus \underline r}$:$$At_u(v,\xi)=c(\xi+\rho(u)v).$$ Here $v,\xi\in
M_{X\setminus \underline r}$, and  the operator   $c$ is the Clifford action of
$M_{X\setminus \underline r}$ on $\bigwedge M_{X\setminus \underline r}$.

Denote by $p,q$ the projections of $M_X= M_{X\setminus \underline r} \times M_{X\cap \underline r}$ on the two factors.
 Consider  $M_{X\setminus \underline r} \times M_{X\cap \underline r}$  as a $G\times G$ manifold.
 We use now the general multiplicative formula  Theorem 3.5 of  \cite{At}, applied to  the groups $G\times G$.

 If $\tau$ is a $G$-transversally elliptic symbol on  $F(\underline r)$, then $ \sigma:=
p^*At_u\odot q^*\tau $ is
a $G\otimes G$-transversally elliptic symbol on $M_{X\setminus \underline r} \times F(\underline r)$ with index the product of the two indices on $G\times G$: $index(p^*At_u\odot q^*\tau)(g_1,g_2)=
index(p^*At_u)(g_1)index(q^*\tau)(g_2)$.
 Now, consider $G$ embedded as the diagonal in $G\times G$. As $\tau$ is
 $G/ G_{\underline r}$ transversally elliptic, and $At_u$ is $G(\underline r)$ transversally elliptic,
   the symbol $\sigma$ remains $G$-transversally elliptic. We thus obtain that
$p^*At_u\odot q^*\tau $ is
a $G$-transversally elliptic symbol on $M_{X}$ with index the product of the two indices.
Let us restrict $\sigma$  to  $T_G^*M_X|_{F(\underline r)}$. This space is the product  $T_G^*F(\underline r)\times M_{X\setminus \underline r} '$   and thus,  for each point  in this product, the vector $v=0$ and  the  symbol $p^*At_u$ coincides with  the  Bott symbol ${\bf c}_{M_{X\setminus \underline r} '} $.   In other words,
by definition of the isomorphism $C_{\underline r}$, the restriction $e^*(\sigma)$  to $T_G^*M_X|_{F(\underline r)}$ is $C_{\underline r}(\tau)$.

Thus we see that $C_{\underline r}^{-1}e^*(\sigma)=\tau$, so that $C_{\underline r}^{-1}e^*$  is surjective and the map $ \tau\mapsto
p^*At_u\odot q^*\tau $ is the desired splitting.

Finally  remark that this splitting   corresponds, under the  index isomorphism,
to the splitting of the spaces
$\tilde {\mathcal F}_i(X) $  using convolution by the partial partition functions $\mathcal  P^{F_{\underline r}}_X$.
This follows from the explicit computation of the index  of the pushed  symbols in the multiplicative formula which we have previously recalled.
\end{remark}

\section{Back to partition functions}

Assume $G$ is a torus and let $M_X:=\oplus _{a\in X}L_a$.

Let
${\rm Cone}(X):=\{\sum _{a\in X} t_a a\,|\, t_a\geq 0\}$
be the cone generated by $X$  in $V$.
Assume  ${\rm Cone}(X)$ is a pointed cone so that we can consider the partition function
$\mathcal P_X$. This partition function counts the number of integral points in the partition polytope
$P(\lambda):=\{t_a, t_a\in \mathbb R_{\geq 0},\,|\, \sum_a t_a a=\lambda\}$.

Let $V_{\rm sing}$ be the union of the cones
${\rm Cone}(Y)$ generated by  the sublists $Y$ of $X$
which do not span $V$. A connected component of $V\setminus V_{\rm sing}$ is called a big cell.
Recall \cite{dpv} that the partition function $\mathcal P_X$ coincide on a big cell  $\mathfrak c$ with an element of $DM(X)$ denoted by $\mathcal P_X^{\mathfrak c}$.

Consider the moment map $\mu_X: M_X\to V$  given by $\mu_X(v)=\sum_{a\in X} \|v_a\|^2 a$.
 Then $\mu_X$ is a proper map.
If $v\in {\rm Cone}(X)$ is in a big cell $\mathfrak c$, $v$ is a regular value of $\mu_X$,
 and the manifold
 $P_v:=\mu_X^{-1}(v)$ is a compact closed submanifold of $M_X^f$ . We recall that
 $M(\mathfrak c):=P_v/G$ is a toric manifold which depends only of $\mathfrak c$.

A neighborhood $\mathcal U_v$ of  $P_v$ in $M_X^f$ is isomorphic to the product $V\times P_v$, so that we have
  $ K_G^0(P_v)\sim K_G^s(\mathcal U_v)$.
  Via the open embedding
  $\mathcal U_v\to M_X^f$, we obtain a map
  $m_v: K_G^0(P_v)\to  K_G^s(M_X^f).$

 Let $I_v$  be the trivial bundle over $P_v$. It is easy to see that the element
 $ m_v(I_v)\in K_G^s(M_X^f)$ depends only of the big cell $\mathfrak c$ where $v$ leaves.
 We denote it by $I_{\mathfrak c}$.

Recall the isomorphism $r:K_G^s(M_X^f)\to DM(X)$  given by combining \ref{vanish}  and the index map \ref{first}.
 The following theorem will be proved in a subsequent article.

 \begin{theorem}\label{back}

 \begin{enumerate}
 \item
 When $\mathfrak c$ runs over all big cells contained in ${\rm Cone}(X)$,
 the elements $I_{\mathfrak c}$ generate  $K_G^s(M_X^f)$ as a $R(G)$ module.\\

 \item We have
 $ r(I_\mathfrak c)= \mathcal P_X^{\mathfrak c}.$

 \end{enumerate}
 \end{theorem}

 The proof of the second item of this theorem follows right away from
the free action property of the index when the toric manifold $M(\mathfrak c)$ is smooth.
  Indeed, choose an Hermitian structure on $M_X$.
 The space $T^*_G P_v$ has the structure of an Hermitian vector bundle.
 Its fiber at the point $p\in P_v$ is isomorphic to the ``horizontal tangent space'' $H_p$.
 As in Example \ref{CR}, we can construct the transversally elliptic symbol
 $\sigma_v(p,\xi)=c(\xi^1)$, where $\xi^1$ is the projection of $\xi$ on $H_p$ and $c$ the Clifford action of $H_p$ on  $\Lambda H_p$.
 Via the closed embedding $i_v: P_v \to M_X^f$, we obtain an element  $(i_v)_{!}(\sigma_v)$ in $K^0(T^*_G M_X^f)$.
 Thus the element $r(I_{\mathfrak c})$ is the index multiplicity of $\sigma_v$.

 Assume $M(\mathfrak c)$ smooth. Then
  $P_v\to P_v/G$ is a principal bundle. If   $\lambda\in \Lambda$,
  we obtain a holomorphic  line bundle $\mathcal O(\lambda)=P_v\times_G \mathbb C_{\lambda}$
  over   $M(\mathfrak c)$.
   By the free action
  property of the index,
  $ind_m(\sigma_{v})(\lambda)$ is the  virtual dimension of the space
  $\sum_{i=0}^{\dim M(\mathfrak c)}(-1)^i H^{0,i}(M(\mathfrak c), \mathcal O(\lambda))$.
 This dimension is polynomial in $\lambda$. On the other hand, when $\lambda\in \mathfrak c$,
 $H^{0,i}(M(\mathfrak c), \mathcal O(\lambda))$ vanishes if $i>0$ and $H^{0,0}(M(\mathfrak c), \mathcal O(\lambda))$  is the number of integral points in the partition polytope $P(\lambda)$, which is given by the function
 $\mathcal P_X(\lambda)$.

     When   $M(\mathfrak c)$ is an orbifold,  we
     deduce  Theorem \ref{back} from the general
     cohomological index theorem for transversally elliptic operators.
     We will discuss this point in a subsequent article.


\begin{thebibliography}{99}

\bibitem{At} {Atiyah M.,}
{Elliptic operators and compact groups} , Springer L.N.M., n. 401, 1974.\newline

\bibitem{Ats}  Atiyah, M. F.; Singer, I. M. The index of elliptic operators. I. Ann. of Math. (2) 87 1968 484--530.\newline




\bibitem{Ats2}  Atiyah, M. F.; Singer, I. M. The index of elliptic operators. III. Ann. of Math. (2) 87 1968 546--604. \newline


\bibitem{BV1}
Berline, Nicole; Vergne, Mich\`ele L'indice \'equivariant des op\'erateurs transversalement elliptiques.   Invent. Math. 124 (1996), no. 1-3, 51--101.\newline









\bibitem{dp1}{De Concini C., Procesi C., } {Topics in hyperplane arrangements, polytopes and box--splines,}
forthcoming book
\newline

\bibitem{dpv}{De Concini C., Procesi C., Vergne M.} {Partition function and generalized Dahmen-Micchelli spaces,}
 Preprint ArXiv 0805.2907
\newline



\bibitem{DM}{
Dahmen W., Micchelli C., }{The number of solutions to linear
Diophantine equations and  multivariate splines, }{\it Trans. Amer.
Math. Soc.} {\bf 308} (1988), no. 2, 509--532.\newline


\bibitem{parver} {Paradan, Paul--\'Emile, Vergne, M. }
{Index of transversally elliptic operators} {\it To appear}  ArXiv 08041225
\newline








\bibitem{zas} {Zaslavsky T.,} {
 Facing up to arrangements: face count formulas for partitions  of space  by hyperplanes,}  Mem. Amer. Math. Soc., vol. 1, no.  1554, 1975.
\newline


\end{thebibliography}
\end{document}